\documentclass[final,onefignum,onetabnum,a4paper]{siamart171218}

\usepackage{lipsum}
\usepackage{amsfonts}
\usepackage{graphicx,epstopdf}
\usepackage[caption=false]{subfig}
\usepackage{algorithmic}
\ifpdf
  \DeclareGraphicsExtensions{.eps,.pdf,.png,.jpg}
\else
  \DeclareGraphicsExtensions{.eps}
\fi

\newsiamremark{remark}{Remark}
\newsiamremark{hypothesis}{Hypothesis}
\crefname{hypothesis}{Hypothesis}{Hypotheses}
\newsiamthm{claim}{Claim}

\headers{Positivity Preserving Limiters for DG Discretizations}{J.J.W. van der Vegt, Yinhua Xia and Yan Xu}

\title{Positivity Preserving Limiters for Time-Implicit Higher Order Accurate Discontinuous Galerkin Discretizations
\thanks{Submitted to the editors November 21, 2018. 
}}

\author{J.J.W. van der Vegt\thanks{Department of Applied Mathematics, University of Twente, P.O. Box 217, 7500AE, Enschede, The Netherlands
  (\email{j.j.w.vandervegt@utwente.nl}).\funding{Research of J.J.W. van der Vegt was supported by fellowships from the University of Science and
Technology of China (USTC), while the author was in residence at USTC, Hefei, Anhui 230026, People's Republic of China.}}
\and Yinhua Xia\thanks{School of Mathematical Sciences, University of Science and Technology of China, Hefei, Anhui 230026, People's Republic of China.
  (\email{yhxia@ustc.edu.cn}).\funding{Research of Yinhua Xia was supported by NSFC grants No. 11471306, 11871449, and a grant from the Science \& Technology on Reliability \& Environmental Engineering Laboratory (No. 6142A0502020817). }}
\and Yan Xu\thanks{Corresponding author.  School of Mathematical Sciences, University of Science and Technology of China, Hefei, Anhui 230026, People's Republic of China.
  (\email{yxu@ustc.edu.cn}) \funding{ Research of Yan Xu was supported by NSFC grants No. 11722112, 91630207.}}
  }

\usepackage{amsopn}

\begin{document}

\maketitle
\begin{abstract}
 Currently, nearly all positivity preserving discontinuous Galerkin (DG) discretizations of partial differential equations are coupled with explicit time integration methods. Unfortunately, for many problems this can result in severe time-step restrictions. The techniques used to develop explicit positivity preserving DG discretizations can, however, not easily be combined with implicit time integration methods. In this paper we therefore present a new approach. Using Lagrange multipliers the conditions imposed by the positivity preserving limiters are directly coupled to a DG discretization combined with a Diagonally Implicit Runge-Kutta time integration method. The positivity preserving DG discretization is then reformulated as a Karush-Kuhn-Tucker (KKT) problem, which is frequently encountered  in constrained  optimization.  Since the limiter is only active in areas where positivity must be enforced it does not affect the higher order DG discretization elsewhere. The resulting non-smooth nonlinear algebraic equations have, however, a different structure compared to most constrained optimization problems. We therefore develop  an efficient active set semi-smooth Newton method that is suitable for the KKT formulation of time-implicit positivity preserving DG discretizations. Convergence of this semi-smooth Newton method is proven using a specially designed quasi-directional derivative of the time-implicit positivity preserving DG discretization. The time-implicit positivity preserving DG discretization is demonstrated for several nonlinear scalar conservation laws, which include the advection, Burgers, Allen-Cahn, Barenblatt, and Buckley-Leverett equations.
  \end{abstract}



\section{Introduction}
\label{sec:intro}
The solution of many partial differential equations frequently must satisfy a maximum principle, or more generally, certain variables must obey a lower and/or upper bound. In this paper we will denote  all these cases with positivity preserving. In particular, if the partial differential equations model physical processes then these bounds  are also crucial to obtain a meaningful physical solution. For example, a density, concentration or pressure in fluid flow must be nonnegative, and a probability distribution should be in the range $[0,1]$. A numerical solution should therefore strictly obey the bounds on the exact solution, otherwise the problem can become ill-posed and the solution would be meaningless. Also, the numerical algorithm can easily become unstable and lack robustness if the numerical solution violates these essential bounds.

In recent years, the development of positivity preserving discontinuous Galerkin (DG) finite element methods therefore has been a very active area of research. The standard approach to ensure that the numerical solution satisfies the bounds imposed by the partial differential equations is to use limiters, but this can easily result in loss of accuracy, especially for higher order accurate discretizations.

In a seminal paper  Zhang and Shu \cite{zhang2010maximum} showed how to design maximum principle and positivity preserving higher order accurate DG methods for first order scalar conservation laws. Their algorithm consists of a several important steps: i.) starting from a bounds preserving solution at time $t_n$ ensure that the element average of the solution satisfies the bounds at the next time level  $t_{n+1}$ by selecting a suitable time step in combination with a monotone first order scheme; ii.) limit the higher order accurate polynomial solution at the quadrature points in each element without destroying the higher order accuracy; iii.) higher order accuracy in time can then be easily obtained using explicit SSP Runge-Kutta methods \cite{shu1988efficient}. This algorithm has been subsequently extended in many directions, e.g. various element shapes, convection-diffusion equation, Euler and Navier-Stokes equations and relativistic hydrodynamics \cite{zhang2012maximum,zhang2013maximum,zhang2010positivity,zhang2011positivity,zhang2017positivity,qin2016bound}. Other approaches to obtain higher order positivity preserving DG discretizations can be found in e.g. \cite{chen2016third,guo2015positivity,guo2017bound}.

All these DG discretizations use, however, an explicit time integration method. For many partial differential equations this results in an efficient numerical discretization, where to ensure  stability the time step is  restricted by the Courant-Friedrichs-Lewy (CFL) condition. On locally dense meshes and for higher order partial differential equations, which often have a time step constraint $\triangle t\leq C h^p$, with $p> 1$ and $h$ the mesh size, these time-explicit algorithms can become computationally very costly. The alternative is to resort to implicit time integration methods, but positivity preserving time-implicit DG discretizations are still very much in their infancy.  Meister and Ortleb developed in \cite{meister2014unconditionally} a positivity preserving DG discretization for the shallow water equations using the Patankar technique \cite{patankar1980numerical}.   Qin and Shu \cite{qin2018implicit} extended the framework in \cite{zhang2010maximum,zhang2010positivity} to implicit positivity preserving DG discretizations of conservation laws in combination with an implicit Euler time integration method. An interesting result of the analysis in \cite{qin2018implicit} is that to ensure positivity in the algorithm of Qin and Shu a lower bound on the time step is required. The approaches in \cite{meister2014unconditionally,qin2018implicit} require, however, a detailed analysis of the time-implicit DG discretization to ensure that the bounds are satisfied and are not so easy to extend to other classes of problems. 

In this paper, we will present  a very different approach to develop positivity preserving higher order accurate DG discretizations that are combined with a Diagonally Implicit Runge-Kutta (DIRK) time integration method. In analogy with obstacle problems we consider the bounds imposed by a maximum principle or positivity constraint as a restriction on the DG solution space. The constraints are then imposed using a limiter and directly coupled to the time-implicit higher order accurate DG discretization using Lagrange multipliers. { The resulting equations are the well-known Karush-Kuhn-Tucker (KKT) equations, which are frequently encountered  in constrained  optimization and solved with a semi-smooth Newton method \cite{facchinei2007finite,ito2008lagrange}, and also used in constrained optimization based discretizations of partial differential equation in e.g. \cite{Bochev2016CMA,Bochev2016bc,evans2009,BochevSJSC2017}.} The key benefit of the approach discussed in this paper, which we denote   KKT Limiter and so far has not been applied to positivity preserving time-implicit DG discretizations, is that no detailed analysis is required   to ensure that the DG discretization preserves the bounds for a particular partial differential equation. They are imposed explicitly and  not part of the DG discretization. Also, since the limiter is only active in areas where positivity must be enforced, it does not affect the higher order DG discretization elsewhere since the Lagrange multipliers will  be zero there. The approach discussed in this paper presents a general framework how to couple DG discretizations with limiters and, very importantly, how to efficiently solve the resulting nonlinear algebraic equations. 

The algebraic equations resulting from the KKT formulation of the positivity preserving time-implicit DG discretization are only semi-smooth. This excludes the use of standard Newton methods since they require $C^1$ continuity \cite{deuflhard2011newton}. The obvious choice would be to use one of the many semi-smooth Newton methods available for nonlinear constrained optimization problems \cite{facchinei2007finite,ito2008lagrange}, but the algebraic equations for the positivity preserving time-implicit DG discretization have a different structure than for most constrained optimization problems.
For instance, the conditions to ensure a non-singular Jacobian \cite{facchinei2007finite} for methods based on the Fischer-Burmeister or related complementarity functions \cite{munson2001semismooth,chen2010smoothing}  are not met by the KKT-limiter in combination with a time-implicit DG discretization.  This frequently results in nearly singular Jacobian matrices, poor convergence and lack of robustness.  We therefore developed an efficient active set semi-smooth Newton method that is suitable for the KKT formulation of time-implicit positivity preserving DG discretizations. Convergence of this semi-smooth Newton method can be proven using a specially designed quasi-directional derivative as outlined in \cite{han1992globally}, see also \cite{ito2008lagrange,ito2009semi}.

The organization of this paper is as follows. In Section \ref{sec:KKT_Limiter} we formulate the KKT equations, followed in Section \ref{sec:semismooth_Newton} by a discussion of an active set semi-smooth Newton method that is suitable to solve the nonlinear algebraic equations resulting from the positivity preserving time-implicit DG discretization. Special attention will be given to the quasi-directional derivative, which is an essential part to ensure convergence of the semi-smooth Newton method. In Section \ref{KKTDGdiscretization} we discuss the DG discretization in combination with an DIRK time integration method and  positivity constraints. In Section \ref{Numerical Experiments} numerical experiments for the advection, Burgers, Allen-Cahn, Barenblatt, and Buckley-Leverett equations are provided. Conclusions are drawn in Section \ref{sec:conclusions}. In the Appendix more details on the quasi-directional derivative are given.

\section{Karush-Kuhn-Tucker limiting approach}
\label{sec:KKT_Limiter}
In this section we will directly couple the bounds preserving limiter to the time-implicit discontinuous Galerkin discretization using Lagrange multipliers. We will denote this approach as the KKT-Limiter.

Define the set
\begin{equation*}
K:=\{x\in\mathbb{R}^n\;\vert\; h(x)=0,\; g(x)\leq 0 \},
\end{equation*}
where $h:\mathbb{R}^n\rightarrow\mathbb{R}^l$ and $g:\mathbb{R}^n\rightarrow\mathbb{R}^m$ are twice continuously differentiable functions denoting, respectively, the $l$ equality and $m$ inequality constraints to be imposed on the DG discretization. The variable $x$ denotes the degrees of freedom and $n$ the number of degrees of freedom in the unlimited DG discretization. For the continuously differentiable function $L:\mathbb{R}^n\rightarrow\mathbb{R}^n$, representing the unlimited discontinuous Galerkin discretization, the KKT-equations are
\begin{subequations}
\begin{align}
\mathcal{L}(x,\mu,\lambda):=L(x)+\nabla h(x)^T\mu+\nabla g(x)^T\lambda&=0,\\
-h(x)&=0,\label{equality_constraint}\\
0\geq g(x)\perp\lambda&\geq 0,\label{compatibility_eq}
\end{align}\label{KKT_equations}
\end{subequations}
with $\mu\in\mathbb{R}^l$, $\lambda\in\mathbb{R}^m$ the Lagrange multipliers.
The compatibility condition \cref{compatibility_eq}  is component-wise equal to:
 \begin{equation*}
 0\geq g_j(x),\qquad\lambda_j\geq 0\quad\text{and}\quad g_j(x)\lambda_j=0,\quad j=1,\cdots,m,
 \end{equation*}
 which is equivalent with
 \begin{equation*}
 \min(-g(x),\lambda)=0,
 \end{equation*}
where the $\min$-function is applied component-wise.
The KKT-equations, with $F(z)\in\mathbb{R}^{n+l+m}$, can now be formulated as
\begin{equation}
0=F(z):=\left(\begin{matrix}
\mathcal{L}(x,\mu,\lambda)\\[5pt]
-h(x)\\[5pt]
\min(-g(x),\lambda)
\end{matrix}\right),\label{KKTnonlineareq}
\end{equation}
where $z:=(x,\mu,\lambda)$. In the next section we will discuss a global active set semi-smooth Newton suitable for the efficient solution of \cref{KKTnonlineareq} in combination with a DIRK-DG discretization. In Section \ref{KKTDGdiscretization} the DG discretization and KKT-Limiter will be presented for a number of scalar conservation laws.

\section{Semi-Smooth Newton Method}\label{sec:semismooth_Newton}
Standard Newton methods assume that $F(z)$ is continuously differentiable \cite{deuflhard2011newton}, but $F(z)$ given by \cref{KKTnonlineareq} is only semi-smooth \cite{facchinei2007finite}. In this section we will present a robust active set semi-smooth Newton method for \cref{KKTnonlineareq} that is suitable for the efficient solution of the KKT-equations resulting from a higher order DG discretization combined with positivity preserving limiters and a Diagonally Implicit Runge-Kutta time integration method \cite{hairer2010solving}.

\subsection{Differentiability concepts}
For the definition of the semi-smooth Newton method we need several more general definitions of derivatives, which will be discussed in this section. For more details, we refer to e.g. \cite{clarke1990optimization,facchinei2007finite,ito2008lagrange,shapiro1990concepts}. Since we use the semi-smooth Newton method directly on the algebraic equations of the limited DIRK-DG discretization we only consider finite dimensional spaces here.

Let  $D\subseteq \mathbb{R}^m$ be an open subset in $\mathbb{R}^m$.
Given $d\in \mathbb{R}^m$, the {directional derivative} of $F:D\rightarrow\mathbb{R}^n$ at $x\in D$ in the direction $d$ is defined as
\begin{equation}
F^\prime(x;d):=\lim_{t\downarrow 0^+}\frac{F(x+t d)-F(x)}{t}.\label{eq:directional_derivative}
\end{equation}
A function $F:D \rightarrow\mathbb{R}^n$ is {locally Lipschitz continuous} if for every $x\in D$ there exists a neighborhood $N_x\subseteq D$ and a constant $C_x$, such that
\begin{equation*}
\vert F(y)-F(z)\vert\leq C_x\vert y-z\vert,\qquad\forall y,z\in N_x.
\end{equation*}
If $F$ is locally Lipschitz on $D$ then according to Rademacher's
theorem $F$ is differentiable almost everywhere with derivative $F^\prime(x)$.
The  B-subdifferential $\partial_BF(x)$ of $F(x)$  is then defined as
\begin{equation*}
\partial_BF(x):=\lim_{\bar{x}\rightarrow x,\bar{x}\in D_F}F^\prime(\bar{x}),
\end{equation*}
with $D_F$ the points where $F$ is differentiable, and the generalized derivative in the sense of Clarke is defined as
\begin{equation*}
\partial F(x):=\text{ convex hull of $\partial_BF(x)$}.
\end{equation*}
For example, $F(x)=\vert x\vert$   at $x=0$ has $\partial_BF(0)=\{-1,1\}$ and  $\partial F(0)=[-1,1]$.
A function $F: D\rightarrow\mathbb{R}^n$ is called semi-smooth if \cite{qi1993nonsmooth}
\begin{equation*}
\lim_{V\in\partial F(x+td^\prime), d^\prime\rightarrow d, t\downarrow 0^+} Vd^\prime\quad\text{exists for all $d\in\mathbb{R}^m$}.
\end{equation*}

A function $F:D\rightarrow \mathbb{R}^n$ is  {Bouligand-differentiable} (B-differentiable) at $x\in D$  if it is directionally differentiable at $x$ and
\begin{equation*}
\lim_{d\rightarrow 0}\frac{F(x+d)-F(x)-F^\prime(x;d)}{\vert d\vert}=0.
\end{equation*}
A locally Lipschitz continuous function $F$ is {B-differentiable} at $x$ if and only if it is directionally differentiable at $x$ \cite{shapiro1990concepts}.

Given $d\in\mathbb{R}^m$, the {Clarke generalized directional derivative} of $F: D\rightarrow
\mathbb{R}^n$ at $x\in D$ in the direction of $d$ is defined by \cite{clarke1990optimization}
\begin{equation*}
F^0(x;d):=\lim_{y\rightarrow x}\sup_{t\downarrow 0^+}\frac{F(y+t d)-F(y)}{t}.
\end{equation*}

\subsection{Global active set semi-smooth Newton method}
For the construction of a global semi-smooth Newton method for \cref{KKTnonlineareq} we will use the merit function $\theta(z)=\frac{1}{2}\vert F(z)\vert^2$, with $z=(x,\mu,\lambda)$. The Clarke directional derivative of $\theta$ and $F$ have the following relation.

Let  $F:D\subseteq\mathbb{R}^p\rightarrow\mathbb{R}^p$,  with $D$ an open set and $p=n+l+m$, be a  locally Lipschitz continuous function then the Clarke generalized directional derivative of $\theta(z)$ can be expressed as \cite{ito2008lagrange}
\begin{equation}
\theta^0(z;d)=\limsup_{y\rightarrow z, t\downarrow 0^+}\frac{(F(z),(F(y+td)-F(y))}{t},\label{thetaclark}
\end{equation}
and there exists an  $F^0:D\times\mathbb{R}^p\rightarrow\mathbb{R}^p$ such that
\begin{equation}
\theta^0(z;d)=(F(z),F^0(z;d))\qquad\text{for $(z,d)\in D\times\mathbb{R}^p$}.\label{thetader1}
\end{equation}
Here $(\cdot,\cdot)$ denotes the Euclidian inner product.
The crucial point in designing a Newton method is to obtain  proper descent directions for the Newton iterations. A possible choice is to use the Clarke derivative $\partial F$ as generalized Jacobian \cite{facchinei2007finite,ito2008lagrange}, but this derivative is in general difficult to compute. In \cite{pang1990newton,pang1991b} it was proposed to use $d$ as the solution of
\begin{equation}
F(z)+F^\prime(z;d)=0,\label{direcderivsearch}
\end{equation}
which  for the KKT-equations results in a mixed linear complementarity problem \cite{pang1991b}. Unfortunately, \cref{direcderivsearch} does not always have a solution, unless additional conditions are imposed. A better alternative is to use the quasi-directional derivative $G$ of $F$
\cite{han1992globally,ito2008lagrange,ito2009semi}.

Let $F:D\subseteq\mathbb{R}^p\rightarrow\mathbb{R}^p$ be directionally differentiable and locally Lipschitz continuous. Assume that $S=\{z\in D\;\vert\; \vert F(z)\vert\leq \vert F(z^0)\vert\}$ is bounded. Then $G:S\times\mathbb{R}^p\rightarrow\mathbb{R}^p$ is called the quasi-directional derivative of $F$ on $S\subset\mathbb{R}^p$ if for all $z,\bar{z}\in S$ the following conditions hold \cite{han1992globally,ito2008lagrange,ito2009semi}
\begin{subequations}
\begin{align}
&(F(z),F^\prime(z;d))\leq (F(z),G(z;d)),\label{quasi_der_cond1}\\[5pt]
&G(z;td)=tG(z;d)\quad\text{for all $d\in\mathbb{R}^p,z\in S$ and $t\geq 0$},\label{quasi_der_cond2}\\[5pt]
&(F(\bar{z}),F^0(\bar{z};\bar{d}))\leq\limsup_{z\rightarrow \bar{z},d\rightarrow \bar{d}}(F(z),G(z;d))\quad\text{for all $z\rightarrow \bar{z},d\rightarrow \bar{d}$}.\label{Gupperbound}
\end{align}\label{quasidirectional_derivative}
\end{subequations}
The search direction $d$  in the semi-smooth Newton method is now the solution of
\begin{equation}
F(z)+G(z;d)=0,\quad\text{with}\; z\in S,d\in\mathbb{R}^p,\label{searcheq}
\end{equation}
which results for the KKT-equations \cref{KKTnonlineareq} in a mixed linear complementarity problem.
Using \cref{thetader1}, \cref{Gupperbound} and \cref{searcheq} this immediately results in the bound
\begin{equation*}
\theta^0(\bar{z};\bar{d})\leq \limsup_{z\rightarrow \bar{z},d\rightarrow \bar{d}}(F(z),G(z;d))=-\lim_{z\rightarrow \bar{z}}\vert F(z)\vert^2=-2\theta(\bar{z}).
\end{equation*}
Hence the search direction $d$ obtained from \cref{searcheq} always provides a descent direction for the merit function $\theta(z)$. The merit function $\theta(z)$ and the quasi-directional derivative $G(z,d)$ can therefore be used to define a global line search semi-smooth Newton algorithm, which is stated in Algorithm \ref{Newton_algorithm}. The key benefit of using the quasi-directional derivative $G$ in this Newton algorithm is that, under the additional assumption $\Vert G(z;d)\Vert\geq L\Vert d\Vert$, with $L>0$ constant, we immediately obtain a proof of the convergence of this algorithm, given by \cite{han1992globally}, Theorem 1.

In the next section we will present the quasi-directional derivative $G$ for the KKT-equations \cref{KKTnonlineareq} and define the active sets used to solve \cref{searcheq} with the semi-smooth Newton algorithm presented in Section \ref{active_set_algorithm}. In Section \ref{KKTDGdiscretization}  Algorithm \ref{Newton_algorithm} will then be used to solve the nonlinear equations resulting from the DG discretization using a KKT-limiter in combination with a Diagonally Implicit Runge-Kutta (DIRK) method.

\subsection{Quasi-directional derivative}
In order to compute the quasi-directional derivative $G$, satisfying the conditions stated in \cref{quasidirectional_derivative}, we first need to compute the directional and  Clarke generalized directional derivatives of the function $F(z)$ defined in \cref{KKTnonlineareq}.

Define $z\in\mathbb{R}^p$, with $p=n+l+m$ as $z=(x,\mu,\lambda)$ with $x\in\mathbb{R}^n$, $\mu\in\mathbb{R}^l$, $\lambda\in\mathbb{R}^m$.
Define $d\in\mathbb{R}^p$ as $d=(u,v,w)$ with $u\in\mathbb{R}^n$, $v\in\mathbb{R}^l$, $w\in\mathbb{R}^m$.
The {directional derivative} $F^\prime(z;d)\in\mathbb{R}^p\times \mathbb{R}^p$  of  $F(z)$ defined in \cref{KKTnonlineareq} in the direction $d$ is equal to
\begin{subequations}
\begin{alignat}{2}
F_i^\prime(z;d)&=D_x\mathcal{L}_i(z)\cdot u+D_\mu \mathcal{L}_i(z)\cdot v+D_\lambda \mathcal{L}_i(z)\cdot w,&\qquad &i\in N_n,\\
F_{i+n}^\prime(z;d)&=-D_xh_i(x)\cdot u,&& i\in N_l,\\
F_{i+n+l}^\prime(z;d)&=-D_xg_i(x)\cdot u,&& i\in \alpha(z),\\
&=\min(-D_xg_i(x)\cdot u,w_i),&& i\in\beta(z),\label{directionaldermin}\\
&=w_i,&& i\in\gamma(z),
\end{alignat}\label{directionalderivative}
\end{subequations}
where the following sets are used
\begin{align*}
N_q&=\left\{j\in\mathbb{N}\;\vert\;1\leq j\leq q\right\},\\
\alpha(z)&=\left\{j\in\mathbb{N}_m\;\vert\;\lambda_j>-g_j(x)\right\},\\
\beta(z)&=\left\{j\in\mathbb{N}_m\;\vert\;\lambda_j=-g_j(x)\right\},\\
\gamma(z)&=\left\{j\in\mathbb{N}_m\;\vert\;\lambda_j<-g_j(x)\right\},
\end{align*}
with $q=n$ or $q=l$. The calculation of most of the terms in \cref{directionalderivative} is straightforward, except \cref{directionaldermin}, which can be computed using a Taylor series expansion of the arguments of $\min(-g_i(x),\lambda_i)$ in the limit  of the directional derivative \cref{eq:directional_derivative}, combined with  the relation $\min(a+b,a+d)-\min(a,a)=\min(b,d)$ and the fact that $i\in\beta(z)$.

The Clarke Generalized derivative of $F(z)$ can be computed using the relations \cref{thetaclark}-\cref{thetader1} and is equal to
\begin{subequations}
\begin{alignat}{2}
F_i^0(z;d)&=D_x\mathcal{L}_i(z)\cdot u+D_\mu \mathcal{L}_i(z)\cdot v+D_\lambda \mathcal{L}_i(z)\cdot w,&\quad &i\in N_n,\\
F_{i+n}^0(z;d)&=-D_xh_i(x)\cdot u,&&\hspace*{-30pt} i\in N_l,\\
F_{i+n+l}^0(z;d)&=-D_xg_i(x)\cdot u,&&\hspace*{-30pt} i\in \alpha(z),\label{Clarke_direc_deriva}\\
&=\max(-D_xg_i(x)\cdot u,w_i),&&\hspace*{-30pt} i\in\beta(z),F_{i+n+l}(z)>0,\label{Clarke_direc_derivb}\\
&=\min(-D_xg_i(x)\cdot u,w_i),&&\hspace*{-30pt} i\in\beta(z),F_{i+n+l}(z)\leq 0,\label{Clarke_direc_derivc}\\
&=w_i,&&\hspace*{-30pt} i\in\gamma(z).
\end{alignat}\label{Clarke_direc_deriv}
\end{subequations}
The calculation of \cref{Clarke_direc_derivb,Clarke_direc_derivc} in $F^0(z;d)$ is non-trivial and is detailed in Appendix \ref{Appendix1}.

Using the results for the directional derivative and  the Clarke generalized directional derivative we can now state a quasi-directional derivative $G:D\times\mathbb{R}^p\rightarrow\mathbb{R}^p$, satisfying the conditions \cref{quasidirectional_derivative}, which   for any $\delta>0$ is equal to
\begin{subequations}
\begin{alignat}{2}
G_i(z;d)&=D_x\mathcal{L}_i(z)\cdot u+D_\mu \mathcal{L}_i(z)\cdot v+D_\lambda \mathcal{L}_i(z)\cdot w,&\quad &i\in N_n,\\
G_{i+n}(z;d)&=-D_xh_i(x)\cdot u,&&\hspace*{-40pt} i\in N_l,\\
G_{i+n+l}(z;d)&=-D_xg_i(x)\cdot u,&&\hspace*{-40pt} i\in \alpha_\delta(z),\\
&=\max(-D_xg_i(x)\cdot u,w_i),&&\hspace*{-40pt} i\in\beta_\delta(z),F_{i+n+l}(z)>0,\\
&=\min(-D_xg_i(x)\cdot u,w_i),&&\hspace*{-40pt} i\in\beta_\delta(z),F_{i+n+l}(z)\leq 0,\\
&=w_i,&&\hspace*{-40pt} i\in\gamma_\delta(z),
\end{alignat}\label{quasidirectionDeriv1}
\end{subequations}
with the sets
\begin{align*}
\alpha_\delta(z)&=\left\{j\in\mathbb{N}_m\;\vert\;\lambda_j>-g_j(x)+\delta\right\},\\
\beta_\delta(z)&=\left\{j\in\mathbb{N}_m\;\vert\;-g_j(x)-\delta\leq\lambda_j\leq -g_j(x)+\delta\right\},\\
\gamma_\delta(z)&=\left\{j\in\mathbb{N}_m\;\vert\;\lambda_j<-g_j(x)-\delta\right\}.
\end{align*}
The main benefit of introducing the $\delta$-dependent sets is that  in practice it is hard to test for the set $\beta(z)$, which  would generally be ignored in real computations due to rounding errors. One would then miss a number of important components in the quasi-directional derivative, which can significantly affect the performance of the Newton algorithm. The set $\beta_\delta$ gives, however, a computational well defined quasi-directional derivative $G(z;d)$. In Appendix \ref{Appendix2} a proof is given that $G(z;d)$ satisfies the conditions stated in \cref{quasidirectional_derivative}, which is the condition required in  \cite{han1992globally}, Theorem 1, to ensure convergence of the semi-smooth Newton method.

The formulation of the quasi-directional derivative $G$ \cref{quasidirectionDeriv1} is, however, not directly useful as a Jacobian in the semi-smooth Newton method due to the $\max$ and $\min$ functions. In order to eliminate these functions we  introduce the following sets
\begin{align*}
I_{\beta_\delta}^{11}(z,d)&:=\{i\in\beta_\delta(z)\;\vert\; F_{i+n+l}(z)>0,-D_xg_i(x)\cdot u>w_i\},\\
I_{\beta_\delta}^{12}(z,d)&:=\{i\in\beta_\delta(z)\;\vert\; F_{i+n+l}(z)>0,-D_xg_i(x)\cdot u\leq w_i\},\\
I_{\beta_\delta}^{21}(z,d)&:=\{i\in\beta_\delta(z)\;\vert\; F_{i+n+l}(z)\leq 0,-D_xg_i(x)\cdot u>w_i\},\\
I_{\beta_\delta}^{22}(z,d)&:=\{i\in\beta_\delta(z)\;\vert\; F_{i+n+l}(z)\leq 0,-D_xg_i(x)\cdot u\leq w_i\},
\end{align*}
and define
\begin{subequations}
\begin{align}
I_\delta^1(z,d)&:=\alpha_\delta(z)\cup I_{\beta_\delta}^{11}(z,d)\cup I_{\beta_\delta}^{22}(z,d),\\
I_\delta^2(z,d)&:=\gamma_\delta(z)\cup I_{\beta_\delta}^{12}(z,d)\cup I_{\beta_\delta}^{21}(z,d).
\end{align}\label{active_sets}
\end{subequations}
The quasi-directional derivative $G(z;d)$ can now be written in a  form suitable to serve as a Jacobian in the active set semi-smooth Newton method defined in Algorithm \ref{Newton_algorithm} to solve \cref{KKTnonlineareq}
\begin{equation*}
G(z;d)=\widehat{G}(z)d,
\end{equation*}
with
\begin{equation}
\widehat{G}(z)=\left(\begin{matrix}
D_x\mathcal{L}_i(z)\vert_{i\in N_n} & D_\mu \mathcal{L}_i(z)\vert_{i\in N_n} & D_\lambda \mathcal{L}_i(z)\vert_{i\in N_n} \\
-D_x h_i(x)\vert_{i\in N_l}&0&0\\
-D_x g_i(x)\vert_{i\in I^1_\delta(z,d)} & 0 & \delta_{ij}\vert_{i,j\in I_\delta^2(z,d)}
\end{matrix}\right)\in\mathbb{R}^{p\times p},\label{eq:quasi-direc-deriv}
\end{equation}
with $\delta_{ij}$ the Kronecker symbol. By  updating  the sets $I_\delta^1(z;d)$ and $I_\delta^2(z;d)$ as part of the Newton method the complementary problem \cref{searcheq} is simultaneously solved with the solution of \cref{KKTnonlineareq}. In general, after a few iterations  the proper sets $I^{1,2}_\delta(z;d)$ will be found and the semi-smooth Newton method then converges like a regular Newton method. Also, one should note that {\it only} the contribution $D_x\mathcal{L}_i(z)$ in \cref{eq:quasi-direc-deriv} depends on the DG discretization in $\mathcal{L}_i(z)$. Hence, the KKT-Limiter provides a  general framework to impose limiters on time-implicit numerical discretizations and could for instance also be applied to time-implicit finite volume discretizations.

\subsection{Active set semi-smooth Newton algorithm}\label{active_set_algorithm}
\algsetup{indent=2em}
\begin{algorithm}[tbh]
\small
\caption{Active set semi-smooth Newton method}
\label{Newton_algorithm}
\begin{algorithmic}[1]
\STATE (A.0) ({\it Initialization}) Let $\bar{\alpha}\geq 0$, $\beta,\gamma\in(0,1)$, $\sigma\in(0,\bar{\sigma})$, $\delta>0$ and $b>C\in\mathbb{R}^+$ arbitrarily large, but bounded. Choose $z^0, d^0\in\mathbb{R}^p$ and  tolerance $\epsilon$.
\STATE (A.1) Scale $z^0$.
\STATE (A.2) ({\it Newton method})
\FOR{$k=0,1,\cdots$ until $\Vert F(z^k)\Vert\leq\epsilon$ \AND $\Vert d^k\Vert\leq\epsilon$}
\STATE Compute the quasi-directional derivative matrix $\widehat{G}_k:=\widehat{G}(z^k)$ given by \cref{eq:quasi-direc-deriv} and the active sets $I^1_\delta(z;d)$, $I^2_\delta(z;d)$ of $\widehat{G}_k$ given by \cref{active_sets}.
\STATE Apply row-column scaling  to $ (\widehat{G}^T_k\widehat{G}_k+\bar{\alpha}\Vert F(z^k)/F(z^0)\Vert I)$, with $I$ the identity matrix,  such that the matrix has a norm $\Vert\cdot\Vert_{L^\infty}\cong 1$.
\IF {there exists a solution $h^k$ to
 \begin{equation}
 (\widehat{G}^T_k\widehat{G}_k+\bar{\alpha}\Vert F(z^k)/F(z^0)\Vert I)h^k=-\widehat{G}_k^TF(z^k),\label{leastsq}
 \end{equation}
 with $\vert h^k\vert\leq b\vert F(z^k)\vert$ \AND
 $$\vert F(z^k+h^k)\vert<\gamma\vert F(z^k)\vert,$$}
\STATE Set $d^k=h^k$, $z^{k+1}=z^k+d^k$, $\alpha_k=1$ and $m_k=0$.
\ELSE
\STATE Choose $d^k=h^k$.
\STATE Compute $\alpha_k=\beta^{m_k}$, where $m_k$ is the first positive integer $m$ for which,
 $$\theta(z^k+\beta^{m_k}d^k)-\theta(z^k)\leq-\sigma\beta^m\theta(z^k).$$
 \vspace*{-10pt}
 \STATE Set $z^{k+1}=z^k+\alpha_kd^k$.
\ENDIF
\ENDFOR
\end{algorithmic}\label{algorithm_1}
\end{algorithm}

As default values we use in Algorithm \ref{algorithm_1}  $\bar{\alpha}=10^{-12}$, $\beta=\gamma=\frac{1}{2}$, $\sigma=10^{-9}$, $\delta=10^{-12}$ and $\epsilon=10^{-8}$.

An important aspect of Algorithm \ref{Newton_algorithm} is that we simultaneously solve the mixed linear complementarity equations \cref{searcheq} for the search direction $d$ as part of the global Newton method using an active set technique. This was motivated by \cite{harker1990damped} and will reduce the mixed linear complementarity problem \cref{searcheq} into a set of linear equations. The use of the active set technique is also based on the observation in \cite{ito2009semi} of the close relation between an active set Newton method and a semi-smooth Newton method. After the proper sets $I^1_\delta(z;d)$, $I^2_\delta(z;d)$ are obtained for the quasi-directional derivative $G(z;d)$ the difference with a  Newton method for smooth problems \cite{deuflhard2011newton} will be rather  small. The mixed linear complementarity problem can, however, have one, multiple or no solutions and, in order to deal also with cases where the matrix $G$ is poorly conditioned, we will use a minimum norm least squares or Gauss-Newton method  to solve the algebraic equations \cref{leastsq}.

For the performance of a Newton algorithm proper scaling of the variables is crucial.  Here, we use  the approach outlined in \cite{deuflhard2011newton} and the Newton method is applied directly to the scaled variables. Also, the matrix $\widehat{G}^T_k\widehat{G}_k+\bar{\alpha}\Vert F(z^k)/F(z^0)\Vert I$ in the Newton method will have a much larger condition number  than the matrix $\widehat{G}_k$. In order to improve the conditioning of this matrix  we use simultaneous iterative row and column scaling in the $L^\infty$-matrix norm, as described in \cite{amestoy2008parallel}. This algorithm very efficiently scales the rows and columns such that an $L^\infty$-matrix norm approximately equal to one is obtained. This gives a many orders of magnitude reduction in the matrix condition number and generally  reduces the condition number of the matrix \cref{leastsq} to the same order as the condition  number of the original matrix $\widehat{G}_k$.

\section{KKT-Limiter DG discretization}\label{KKTDGdiscretization}
Given a domain $\Omega\subseteq\mathbb{R}^d$,  $d={\rm dim}(\Omega)$, $d=1,2$, with Lipschitz continuous boundary $\partial\Omega$.
As general model problem we consider the following second order nonlinear scalar equation
\begin{equation}
\frac{\partial u}{\partial t}+\nabla\cdot F(u)+G(u)-\nabla\cdot(\nu(u)\nabla u)=0,\label{eq:2ndconservation_law}
\end{equation}
with $u(x,t):\mathbb{R}^d\times\mathbb{R}^+\rightarrow\mathbb{R}$ a scalar quantity, $F(u):\mathbb{R}\rightarrow\mathbb{R}^d$ the flux, $G(u):\mathbb{R}\rightarrow\mathbb{R}$ a reaction term and $\nu(u):\mathbb{R}\rightarrow\mathbb{R}^+$ a nonlinear diffusion term.
By selecting different functions $F, G$ and $\nu$ in \cref{eq:2ndconservation_law} we will demonstrate in Section \ref{Numerical Experiments}  the KKT-Limiter on various model problems that impose different positivity constraints on the solution.

For the DG discretization we introduce the auxiliary variable $Q\in\mathbb{R}^d$  and rewrite
  \cref{eq:2ndconservation_law} as a first order system of conservation laws
\begin{subequations}
\begin{align}
\frac{\partial u}{\partial t}+\nabla\cdot F(u)+G(u)-\nabla\cdot(\nu(u)Q)&=0,\\
Q-\nabla u&=0.
\end{align}\label{eq:1stconservation_law}
\end{subequations}

\subsection{DG discretization}\label{DGdiscretization}
Let $\mathcal{T}_h$ be a tessellation of the domain $\Omega$ with shape regular line or quadrilateral elements $K$ with maximum diameter $h>0$.  The total number of elements in $\mathcal{T}_h$ is  $N_K$. We denote the union of the set of all  boundary faces $\partial K$, $K\in\mathcal{T}_h$, as $\mathcal{F}_h$, all internal faces  ${\mathcal F}_h^i$ and the boundary faces as ${\mathcal F}^b_h$, hence $\mathcal{F}_h=\mathcal{F}_h^i\cup\mathcal{F}_h^b$. The elements  connected to each side of a face $S\in\mathcal{F}_h$ are denoted by the indices $L$ and $R$, respectively.  { For the  KKT-Limiter it is important to use orthogonal basis functions, see Section \ref{Limiter_Constraints}. In this paper $\mathcal{P}_p(K)$ represent tensor product Legendre polynomials of degree $p$ on $d$-dimensional rectangular elements $K\in\mathcal{T}_h$, when $K$ is mapped to the reference element $(-1,1)^d$. For general elements one can use Jacobi polynomials with  proper weights to obtain an orthogonal basis, see \cite{karniadakis2013spectral}, Section 3.2.}
Next, we define the finite element spaces
\begin{align*}
V_h^p&:=\Big\{v\in L^2(\Omega)\;\vert\; v\vert_K\in \mathcal{P}_p(K),\forall K\in\mathcal{T}_h\Big\},\\
W_h^p&:=\Big\{v\in (L^2(\Omega))^d\;\vert\; v\vert_K\in (\mathcal{P}_p(K))^d,\forall K\in\mathcal{T}_h\Big\},
\end{align*}
with $L^2(\Omega)$ the Sobolev space of square integrable functions.
Equation (\ref{eq:1stconservation_law}) is discretized using the Local Discontinuous Galerkin discretization from  \cite{cockburn1998local}.
Define $L^1_h: V_h^p\times W_h^p\times V_h^p\rightarrow \mathbb{R}$ and $L^2_h: V_h^p\times W_h^p\rightarrow \mathbb{R}$ as
\begin{align}
L^1_h(u_h,Q_h;v):=&-\big(F(u_h)-\nu(u_h)Q_h,\nabla_h v\big)_\Omega+\big(G(u_h),v\big)_\Omega\nonumber\\[4pt]
+&\sum_{S\in\mathcal{F}_h^i}\big(H(u_h^L,u_h^R;n^L)-\widehat{\nu(u_h)}n^L\cdot\widehat{{Q}_h},v^L-v^R\big)_S\nonumber\\
+&\sum_{S\in\mathcal{F}_h^b}\big(H(u_h^L,u_h^b;n^L)-\widehat{\nu(u_h)}n^L\cdot Q_{h}^b,v^L\big)_S,\label{L1_discretization}\\
L^2_h(u_h;w):=&\big(u_h,\nabla_h\cdot w\big)_\Omega
-\sum_{S\in\mathcal{F}_h^i}\big(\widehat{u_h}n^L,w^L-w^R\big)_S\nonumber\\
-&\sum_{S\in\mathcal{F}_h^b}\big(u_h^bn^L,w^L\big)_S,\nonumber
\end{align}
where $(\cdot,\cdot)_D$ is the $L^2(D)$ inner product, $\nabla_h$ the element-wise nabla operator and the superscript $b$ refers to boundary data.
Here $n^L\in\mathbb{R}^d$ is the exterior unit normal vector at the boundary of the element  $L\in\mathcal{T}_h$ that is connected to face $S$.  The numerical flux $H$ is the Lax-Friedrichs flux
 \begin{equation*}
H(u_h^L,u_h^R;n)=\frac{1}{2}\big(n\cdot(F(u_h^L)+F(u_h^R))-C_{LF}(u_h^R-u_h^L)\big),
\end{equation*}
with Lax-Friedrichs coefficient $C_{LF}=\sup_{u_h\in[u_h^L,u_h^R]}\vert\frac{\partial}{\partial u_h}(n\cdot F(u_h))\vert$. For $\widehat{{Q}_h}$ and $\widehat{u_h}$ we use the alternating fluxes
\begin{subequations}
\begin{align}
\widehat{Q_h}&=(1-\alpha)Q_{h}^L+\alpha Q_{h}^R,\\
\widehat{u_h}&=\alpha u_h^L+(1-\alpha)u_h^R,
\end{align}\label{upwind_flux}
\end{subequations}
with $0\leq\alpha\leq 1$. The numerical flux for the nonlinear diffusion is defined as
\begin{equation*}
\widehat{\nu(u_h)}=\frac{1}{2}(\nu(u_h^L)+\nu(u_h^R)).\end{equation*}
For $t\in(0,T]$ the semi-discrete DG formulation for \cref{eq:1stconservation_law} now can be expressed as: Find $u_h(t)\in V_h^p$,  $Q_h(t)\in W_h^p$, such that for all $v\in V_h^p$, $w\in W_h^p$,
\begin{subequations}
\begin{align}
\Big(\frac{\partial u_h}{\partial t},v\Big)_\Omega+L^1_h(u_h,Q_h;v)&=0,\label{eq:uhequation}\\
(Q_h,w)_\Omega+L^2_h(u_h;w)&=0.\label{eq:Qhequation}
\end{align}
\end{subequations}
These equations are discretized in time with a Diagonally Implicit Runge-Kutta (DIRK) method \cite{hairer2010solving}. The main benefit of the DIRK method is that the Runge-Kutta stages can be computed successively, which significantly reduces the computational cost and memory overhead.

{ 
We represent $u_h$ and $Q_h$ in each element $K\in\mathcal{T}_h$, respectively, as $u_h\vert_K=\sum_{j=1}^{N_u}\widehat{U}_j^K\phi_j^K$ and $Q_h\vert_K=\sum_{j=1}^{N_Q}\widehat{Q}_j^K\psi_j^K$, with basis functions $\phi_j^K\in\mathcal{P}_p(K)$, $\psi_j^K\in\big(\mathcal{P}_p(K)\big)^d$ and DG coefficients $\widehat{U}_j^K\in\mathbb{R}$, $\widehat{Q}_j^K\in\mathbb{R}^{d}$. After replacing the test functions $v\in V_h^p$ in \cref{eq:uhequation} and $w\in W_h^p$ \cref{eq:Qhequation} with, respectively, the independent basis functions $\phi_i^K\in\mathcal{P}_p(K)$, $i=1,\cdots,N_u$, and $\psi_i^K\in\big(\mathcal{P}_p(K)\big)^d$, $i=1,\cdots,N_Q$, we obtain the algebraic equations for the DG discretization.

In order to simplify notation we introduce 
$\widehat{L}_h^1(\widehat{U},\widehat{Q})=L_h^1(u_h,Q_h;\phi)\in\mathbb{R}^{N_uN_K}$ and $\widehat{L}_h^2(\widehat{U})=L_h^2(u_h;\psi)\in\mathbb{R}^{dN_QN_K}$, with $N_K$ the number of elements in $\mathcal{T}_h$ and $\phi=\phi_i^K$, $\psi=\psi_i^K$ the basis functions in  element $K$. The algebraic equations for the DIRK stage vector 
$\widehat{K}^{(i)}\in\mathbb{R}^{N_uN_K}$, $i=1,\cdots,s$ with the DG coefficients, then can be expressed as
\begin{align}
\widehat{L}_h(\widehat{K}^{(i)}):=&M_1\big(\widehat{K}^{(i)}-\widehat{U}^n\big)+\triangle t\sum_{j=1}^i a_{ij}\widehat{L}^1_h\big(\widehat{K}^{(j)},-M_2^{-1}\widehat{L}^2_h(\widehat{K}^{(j)})\big)=0.\label{DIRK_stage}
\end{align}
Here we eliminated the DG coefficients for the auxiliary variable $Q_h$ using \cref{eq:Qhequation}.}  The matrices $M_1\in\mathbb{R}^{N_uN_K\times N_uN_K}$, $M_2\in\mathbb{R}^{dN_QN_K\times dN_QN_K}$ are block-diagonal mass matrices since we use orthogonal basis functions and $n$ denotes the index of  time level $t=t_n$.

The coefficients $a_{ij}$ are the coefficients in the Butcher tableau, which determine the properties of the Runge-Kutta method \cite{hairer2010solving}. For DIRK methods $a_{ij}=0$ if $j>i$. The following DIRK methods are used: for basis functions with polynomial order $p=1$ \cite{Alexander1977diagonally}, Page 1012, Theorem 5, first method with $\alpha=1-\frac{1}{2}$; $p=2$ \cite{skvortsov2006diagonally}, Page 2117 (top); $p=3$ \cite{Alexander1977diagonally} Page 1012, Theorem 5, second method, see also \cite{skvortsov2006diagonally}, Page 2117 (top). { The order of accuracy of these DIRK methods is $p+1$ and their coefficients in the Butcher tableau satisfy $a_{sj}=b_j$, $j=1,\cdots,s$, which implies that these methods are stiffly accurate, see \cite{hairer2010solving}, Section IV.6, and the solution of the last DIRK stage is equal to the solution at the new time-step 
\begin{equation*}
\widehat{U}^{n+1}=\widehat{K}^{(s)}.
\end{equation*}
Since each DIRK stage vector must satisfy the positivity constraints this then also immediately applies to the solution at time $t_{n+1}$.

The Jacobian  $D_x\mathcal{L}(\widehat{K}^{(i)})\in\mathbb{R}^{N_uN_K\times N_uN_K}$, with $x=\widehat{K}^{(i)}$, in the quasi-directional derivative $G$ \cref{eq:quasi-direc-deriv} of DIRK stage $i$ of the unlimited DIRK-DG discretization \cref{DIRK_stage}  is now equal to
\begin{equation*}
D_x\mathcal{L}(\widehat{K}^{(i)})=M_1+\triangle t a_{ii}\Big(\frac{\partial L_h^1}{\partial\widehat{K}^{(i)}}-\frac{\partial L_h^1}{\partial\widehat{Q}^{(i)}}M_2^{-1}\frac{\partial L_h^2}{\partial\widehat{K}^{(i)}}\Big).
\end{equation*}

\subsection{Limiter constraints}\label{Limiter_Constraints}
The limiter constraints for the DG discretization can be imposed directly by defining the inequality constraints in the KKT-equations. In each element $K\in\mathcal{T}_h$ we apply for each DIRK-stage $i=1,\cdots,s$, the following inequality constraints
\begin{itemize}
\item[{\it i.}] {\it Positivity constraint}
\begin{equation}
g^K_{1,k}(\widehat{K}^{K,(i)})=u_{\min}-\sum_{q=1}^{N_u}\widehat{K}^{K,{(i)}}_q\phi_q^K(x_k),\quad k=1,\cdots,N_p,\label{positivity_limiter}
\end{equation}
\item[{\it ii.}] {\it Maximum constraint}
\begin{equation}
g^K_{2,k}(\widehat{K}^{K,(i)})=\sum_{q=1}^{N_u}\widehat{K}^{K,{(i)}}_q\phi_q^K(x_k)-u_{\max},\quad k=1,\cdots,N_p.\label{maximum_limiter}
\end{equation}
Here the superscript  $K$ refers to  element $K\in\mathcal{T}_h$, and $(i)$ is the $i$-th DIRK-stage. The points $x_k$, $k=1,\cdots,N_p$, are the points in element $K$ where the inequality constraints are imposed and $u_{\min}$ and $u_{\max}$ denote, respectively, the allowed minimum and maximum value of $u$. The inequality constraints are imposed using the Lagrange multiplier $\lambda$, see \cref{compatibility_eq}.
\smallskip

\item[{\it iii.}] {\it Conservation constraint}
\smallskip

Since the basis functions $\phi_j^K,j=1,\cdots,N_u$ are orthogonal in each element $K$, we have $(1,\phi_j^K)_K=0$, for $j=2,\cdots, N_u$. Hence, at each Runge-Kutta stage $i$, limiting the DG coefficients  $\widehat{K}^{K,(i)}_j$ with $j=2,\cdots,N_u$ has no effect on the element average $\bar{u}^{K,(i)}_h=\frac{1}{\vert K\vert}(u_h^{(i)},1)_K=\widehat{K}^{K,(i)}_1$,  with  $u_h^{(i)}$ the solution at stage $i$, and therefore does not influence the conservation properties of the DG discretization. 
\smallskip

Limiting the DG coefficients $\widehat{K}^{K,(i)}_1$ can, however, effect the conservation properties of the DG discretization since $\bar{u}^{K,(i)}_h=\widehat{K}^{K,(i)}_1$. In order to ensure local conservation we therefore need to impose in each element the local conservation constraint
\begin{align}
h^K\big(\widehat{K}^{K,(i)}\big)&=\widehat{L}_{h,1}^K(\widehat{K}^{(i)})\nonumber\\
&=\vert K\vert \big(\widehat{K}^{K,(i)}_1-\widehat{U}^{n}_1\big)+(G(u_h^{(i)},\phi_1^K)_K\nonumber\\[3pt]
&+\sum_{S\in\mathcal{F}_h^i\cap\partial K}\big(H(u_h^{L,(i)},u_h^{R,(i)};n^L)-\nonumber\\&\hspace*{58pt}\widehat{\nu(u_h)}n^L\cdot
((1-\alpha)Q_h^{L,(i)}+\alpha Q_h^{R,(i)}),\phi^L_1-\phi^R_1\big)_S\nonumber\\[5pt]
&+\sum_{S\in\mathcal{F}_h^b\cap\partial K}\big(H(u_h^{L,(i)},u_h^b;n^L)-\widehat{\nu(u_h)}n^L\cdot Q_{h}^b,\phi^L_1\big)_S,\label{cons_constr}
\end{align}
with $\widehat{L}_{h,1}^K$ the equation for the element mean in element $K$  in \cref{DIRK_stage}. The conservation constraint \cref{cons_constr} is imposed using the Lagrange multiplier $\mu$, see \cref{equality_constraint}.
The conservation constraint explicitly ensures that  at each Runge-Kutta stage the equation for the element mean $\bar{u}^{K,(i)}_h$ is exactly preserved in each element, hence the KKT limiter does not affect the conservation properties of the DG discretization.
\end{itemize}
\smallskip

The remaining Jacobians $D_xh_i(x) \in\mathbb{R}^{N_K\times N_uN_K}$,  $D_xg_i(x) \in\mathbb{R}^{N_pN_K\times N_uN_K}$, and $D_\mu\mathcal{L}_i(z)\in\mathbb{R}^{N_uN_K\times N_K}$,
$D_\lambda\mathcal{L}_i(z)\in\mathbb{R}^{N_uN_K\times N_pN_K}$, with $x=\widehat{K}^{(i)}$,  in the quasi-directional derivative matrix $\widehat{G}$ \cref{eq:quasi-direc-deriv} are now straightforward to calculate. }

It is important to ensure that the initial solution also satisfies the positivity constraints. An $L^2$-projection of the solution will in general not satisfy these constraints for a non-smooth solution. To ensure that the initial solution also satisfies the positivity constraints we apply a constrained projection using the active set semi-smooth Newton method given by Algorithm \ref{algorithm_1}. The only difference is now that instead of \cref{DIRK_stage} we use $L^2$-projection
\begin{equation*}
\widehat{L}_{h i}(\widehat{U}^0)=M^1\widehat{U}^0-(u_0,\phi_i)_\Omega,
\end{equation*}
and combine this with the positivity constraints \cref{positivity_limiter}-\cref{maximum_limiter}. Here,  $u_0$ denotes the initial solution. As initial solution for the constrained projection we use  in Algorithm \ref{algorithm_1} the standard $L^2$-projection without constraints.

The positivity constraints are imposed at all element quadrature points, since only the solution at these quadrature points is used in the DG discretization. In 1D we use Gauss-Lobatto quadrature rules and in 2D product Gauss-Legendre quadrature rules. Since the number of quadrature points in an element is generally larger than the number of degrees of freedom in an element this will result in an over-determined set of algebraic equations and a rank deficit Jacobian matrix if the number of active constraints in an element is larger than the degrees of freedom $N_u$ in element. In order to obtain in Algorithm \ref{algorithm_1}  accurate search directions  $h^k$   we use the Gauss-Newton method given by \cref{leastsq}. This approaches can efficiently deal with the possible rank deficiency of the Jacobian matrix.

In practice it will not be necessary to apply the inequality constraints in all elements and one can significantly reduce the computational cost and memory overhead by excluding those elements for which it is obvious that they will meet the constraints anyway.

\section{Numerical experiments}\label{Numerical Experiments}
In this section we will discuss a number of numerical experiments to demonstrate the performance of the DIRK-DG scheme with the positivity preserving KKT Limiter. All computations were performed using the default values for the coefficients listed for Algorithm \ref{algorithm_1}, except that for the accuracy tests discussed in Section \ref{accuracy_tests} we use $\epsilon=10^{-10}$. The upwind coefficient $\alpha$ in \cref{upwind_flux} is set to $\alpha=1$. In all 1D computations the local conservation constraint is imposed and satisfied with an error less than $10^{-12}$.
\subsection{Accuracy tests}\label{accuracy_tests}
It is important to investigate if the KKT-limiter negatively affects the accuracy of the DG discretization in case the exact solution is smooth, but where also a positivity preserving limiter is required to ensure that the numerical solution stays within the bounds. To investigate this we conduct the same accuracy tests as conducted in Qin and Shu \cite{qin2018implicit}, Section 5.1. Both the linear advection and inviscid  Burgers equation are considered, which are obtained by setting  $F(u)=u$ and $F(u)=\frac{1}{2}u^2$, respectively, and $G(u)=\nu(u)=0$ in \cref{eq:2ndconservation_law}.

\indent{\it Example 5.1} (steady state solution to linear advection equation). We consider
\begin{equation}
u_t+u_x=\sin^4x,\qquad u(x,0)=\sin^2x,\quad u(0,t)=0,\label{linadvec_eq}
\end{equation}
with outflow boundary condition at $x=2\pi$.
The exact solution $u(x,t)$ is positive for all $t>0$, see \cite{qin2018implicit}. As steady state solution we use the solution at $t=500$, when all residuals are approximately $10^{-16}$. During the computations the CFL number is dynamically adjusted between 10 and 89. For the time integration an implicit Euler method is used. In Tables \ref{linadvec_steady_nolimiter} and \ref{linadvec_steady_limiter} the results of the accuracy tests, without and with the KKT-limiter, are shown. The results in Table \ref{linadvec_steady_limiter} show  that the KKT-limiter does not negatively affect the accuracy. For all test cases the optimal accuracy in the $L^2$- and $L^\infty$-norms is obtained. Also, the limiter is necessary, as can be seen from Table \ref{linadvec_steady_nolimiter}, and  preserves the imposed positivity bound $u_{h\min}=10^{-14}$ for the numerical solution.
\begin{table}[htb]
\centering
\caption{Error table for steady state linear advection equation \cref{linadvec_eq} without limiter.}
\begin{tabular}{c|c|c|c|c|c|c}
  \hline\hline
$p$ & $N$ &$L^2 $ error & Order &$ L^\infty$ error & Order &$\min u_h$ \\\hline
&20      &1.461068e-02&-&2.044253e-02&&-5.169578e-03 \\
&40      &3.702581e-03&1.98&5.287628e-03&1.95&-2.883487e-04 \\
1&80    &9.288342e-04&2.00&1.331962e-03&1.99&-1.208793e-05 \\
&160   &2.324090e-04&2.00&3.336614e-04&2.00&-4.036603e-07  \\
&320    &5.811478e-05&2.00&8.345620e-05&2.00&-1.282064e-08 \\\hline
&20       &9.287703e-04&-&1.776878e-03&-&-4.952018e-05  \\
&40       &1.177042e-04&2.98&2.489488e-04&2.84&-1.627459e-06 \\
2&80    &1.476405e-05&3.00&3.200035e-05&2.96&-5.149990e-08 \\
&160   &1.847107e-06&3.00&4.027944e-06&2.99&-1.614420e-09   \\
&320    &2.309385e-07&3.00&5.043677e-07&3.00&-5.049013e-11 \\\hline
&20      &5.653820e-05&-&1.230308e-04&-&-3.877467e-05  \\
&40    &3.583918e-06&3.98&7.803741e-06&3.98&-1.326415e-06 \\
3&80   &2.247890e-07&3.99&4.950122e-07&3.98&-4.237972e-08  \\
&160   &1.406175e-08&4.00&3.090593e-08&4.00&-1.331692e-09  \\
&320   &8.790539e-10&4.00&1.935324e-09&4.00&-4.167274e-11  \\\hline\hline
\end{tabular}
\label{linadvec_steady_nolimiter}
\end{table}

\begin{table}[htb]
\centering
\caption{Error table for steady state linear advection equation \cref{linadvec_eq} with limiter.}
\begin{tabular}{c|c|c|c|c|c|c}
  \hline\hline
$p$ & $N$ &$L^2 $ error & Order &$ L^\infty$ error & Order &$\min u_h$ \\\hline
&20       &1.464990e-02&-&2.044253e-02&-&9.998946e-15\\
&40      &3.702367e-03&1.98&5.287628e-03&1.95&9.999813e-15\\
1&80    &9.288338e-04&2.00&1.331962e-03&1.99&1.000000e-14\\
&160   &2.324090e-04&2.00&3.336614e-04&2.00&1.000000e-14  \\
&320    &5.811478e-05&2.00&8.345620e-05&2.00&1.000000e-14\\\hline
&20       &9.290268e-04&-&1.776878e-03&-&1.000000e-14\\
&40  &1.177053e-04&2.98&2.489488e-04&2.84&1.000000e-14\\
2&80    &1.476406e-05&3.00&3.200035e-05&2.96&1.000000e-14\\
&160   &1.847107e-06&3.00&4.027944e-06&2.99&1.000000e-14  \\
&320    &2.309385e-07&3.00&5.043677e-07&3.00&1.000000e-14\\\hline
&20       &5.742649e-05&-&1.230309e-04&-&9.999990e-15\\
&40    &3.592170e-06&4.00&7.803745e-06&3.98&1.000000e-14\\
3&80    &2.248562e-07&4.00&4.950122e-07&3.98&1.000000e-14\\
&160   &1.406228e-08&4.00&3.090593e-08&4.00&1.000000e-14  \\
&320    &8.790580e-10&4.00&1.935323e-09&4.00&1.000000e-14\\\hline\hline
\end{tabular}
\label{linadvec_steady_limiter}
\end{table}

\indent{\it Example 5.2} (steady state solution to inviscid Burger's equation). We consider the inviscid Burgers equation
\begin{equation}
u_t+(\frac{1}{2}u^2)_x=\sin^3\Big(\frac{x}{4}\Big),\qquad u(x,0)=\sin^2\Big(\frac{x}{4}\Big),\quad u(0,t)=0,\label{Burgers_eq}
\end{equation}
with outflow boundary condition at $x=2\pi$. The exact solution $u(x,t)$ is positive for all $t>0$, see \cite{qin2018implicit}. As steady state solution we use the solution at $t=20.000$, when all residuals are approximately $10^{-16}$. During the computations the CFL number is dynamically adjusted between 10 and 954. For the time integration an implicit Euler method is used. In Tables \ref{Burgers_steady_nolimiter} and \ref{Burgers_steady_limiter} the results of the accuracy tests, without and with the KKT-limiter, show that the KKT-limiter does not negatively affect the accuracy. For all test cases  optimal accuracy in the $L^2$- and $L^\infty$-norms is obtained. Also, the limiter is necessary and  preserves the imposed positivity bound $u_{h\min}=10^{-14}$ for the numerical solution.
\begin{table}[htb]
\centering
\caption{Error table for steady state inviscid Burgers equation \cref{Burgers_eq} without limiter.}
\begin{tabular}{c|c|c|c|c|c|c}
  \hline\hline
$p$ & $N$ &$L^2 $ error & Order &$ L^\infty$ error & Order &$\min u_h$ \\\hline
&20     &2.110016e-03&-&3.387013e-03&-&-2.347303e-03  \\
&40     &5.230241e-04&2.01&8.577912e-04&1.98&-5.865522e-04 \\
1&80   &1.297377e-04&2.01&2.151386e-04&2.00&-1.466204e-04 \\\hline
&20     &2.122765e-05&-&3.024868e-05&-&-1.048636e-05  \\
&40     &2.623666e-06&3.02&3.731754e-06&3.02&-6.681764e-07 \\
2&80    &3.266401e-07&3.01&4.634046e-07&3.01&-4.196975e-08\\ \hline
&20      &2.985321e-07&-&1.895437e-06&-&1.895437e-06 \\
&40      &1.452601e-08&4.36&1.196963e-07&3.99&1.196963e-07\\
3&80    &7.368455e-10&4.30&7.500564e-09&4.00&7.500564e-09\\
&160    &3.948207e-11&4.22&4.346084e-10&4.11& 4.346084e-10\\\hline\hline
\end{tabular}
\label{Burgers_steady_nolimiter}
\end{table}

\begin{table}[htb]
\centering
\caption{Error table for steady state inviscid Burgers equation \cref{Burgers_eq} with limiter.}
\begin{tabular}{c|c|c|c|c|c|c}
  \hline\hline
$p$ & $N$ &$L^2 $ error & Order &$ L^\infty$ error & Order &$\min u_h$ \\\hline
&20      &2.208009e-03&-&3.637762e-03&-&9.999813e-15 \\
&40     &5.358952e-04&2.04&9.282398e-04&1.97&1.000003e-14\\
1&80    &1.313948e-04&2.03&2.339566e-04&1.99&1.000003e-14\\\hline
&20       &2.116746e-05&-&3.024864e-05&-&1.000003e-14\\
&40      &2.622584e-06&3.01&3.731752e-06&3.02&1.000139e-14\\
2&80    &3.266221e-07&3.01&4.634046e-07&3.01&1.000040e-14\\\hline
&20       &2.985321e-07&-&1.895437e-06&-&1.895437e-06\\
&40      &1.452601e-08&4.36&1.196963e-07&3.99&1.196963e-07\\
3&80    &5.610147e-10&4.70&1.574760e-09&6.25&1.000105e-14\\
&160    &3.232240e-11&4.11& 9.038604e-11&4.12&1.000017e-14\\\hline\hline
\end{tabular}
\label{Burgers_steady_limiter}
\end{table}

\subsection{Time dependent tests}
In this section we will present results of simulations of the linear advection, Allen-Cahn, Barenblatt and Buckley-Leverett equations. The order of accuracy of the DIRK time integration method is always $p+1$, with $p$ the polynomial order of the spatial discretization. The minimum value of the residual $F(z)$ and Newton update $d$ in Algorithm \ref{algorithm_1} to stop the Newton iterations is $\epsilon=10^{-8}$ for each DIRK stage. This is a quite strong stopping criteria and in practice the values are often smaller at the end of each DIRK-stage. It is also important to make sure that the Newton stopping criterion is in balance with the accuracy required for the constraints. If the algebraic equations are not solved sufficiently accurate then it is not likely that the KKT-constraints will be satisfied.

The time step for the DIRK method is dynamically computed, based on the CFL or diffusion number. If the Newton method does not converge within a predefined number of iterations, then the computation for the time step will be restarted with $\triangle t/2$. This is generally more efficient than conducting many Newton iterations. In the next time step the time step  will then be increased to $1.2\triangle t$, until the maximum CFL-number is obtained. In practice, depending on the severity of the nonlinearity, the time step will be constantly adjusted during the computations.

\begin{figure}[tbhp]
  \centering
  \hspace*{-5pt}\subfloat[$p=1$]{\label{Linear_Advection1}\includegraphics[width=0.48\textwidth,height=0.25\textheight]{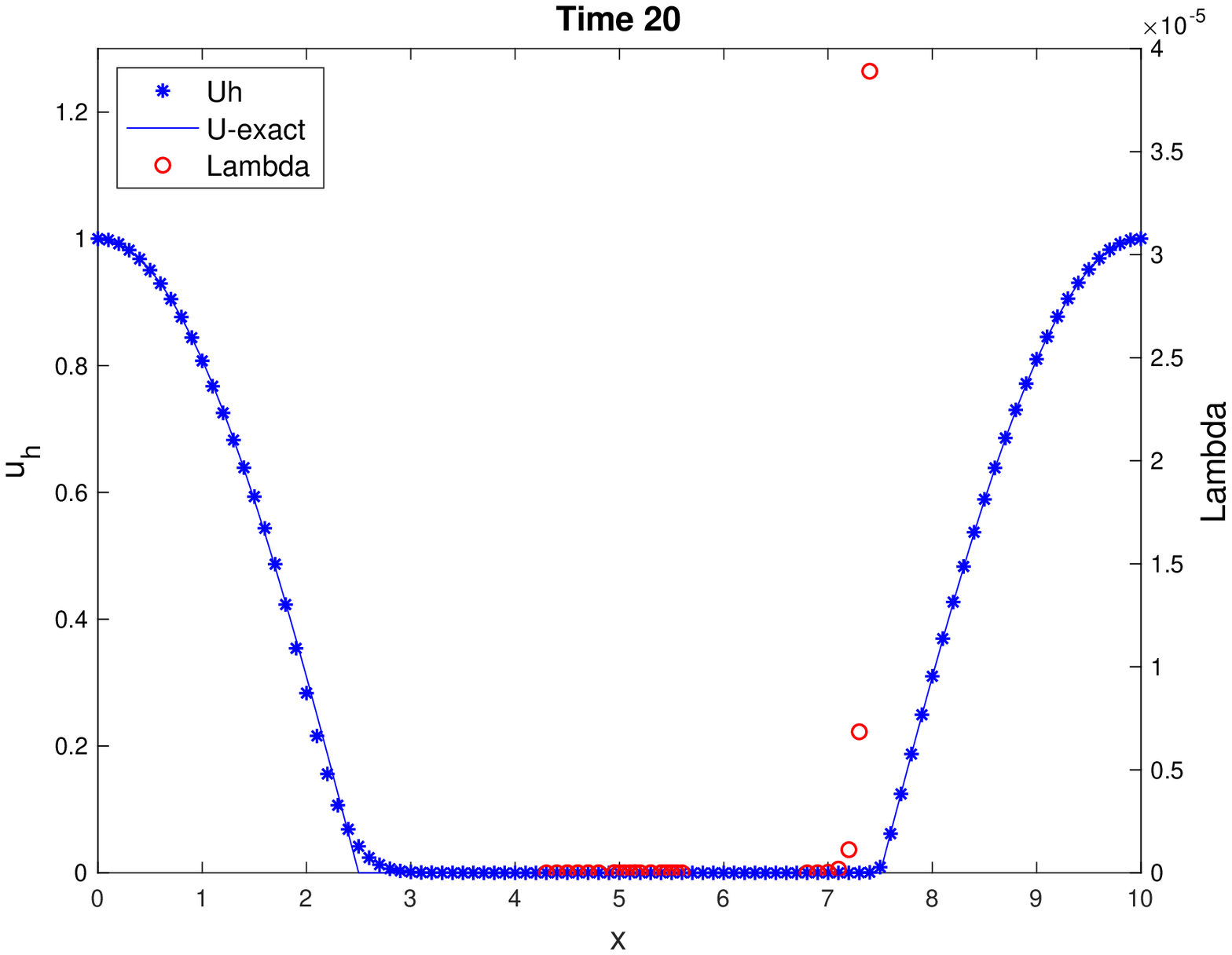}}\hspace*{4pt}
  \subfloat[$p=1$]{\label{Linear_Advection2}\includegraphics[width=0.45\textwidth,height=0.25\textheight]{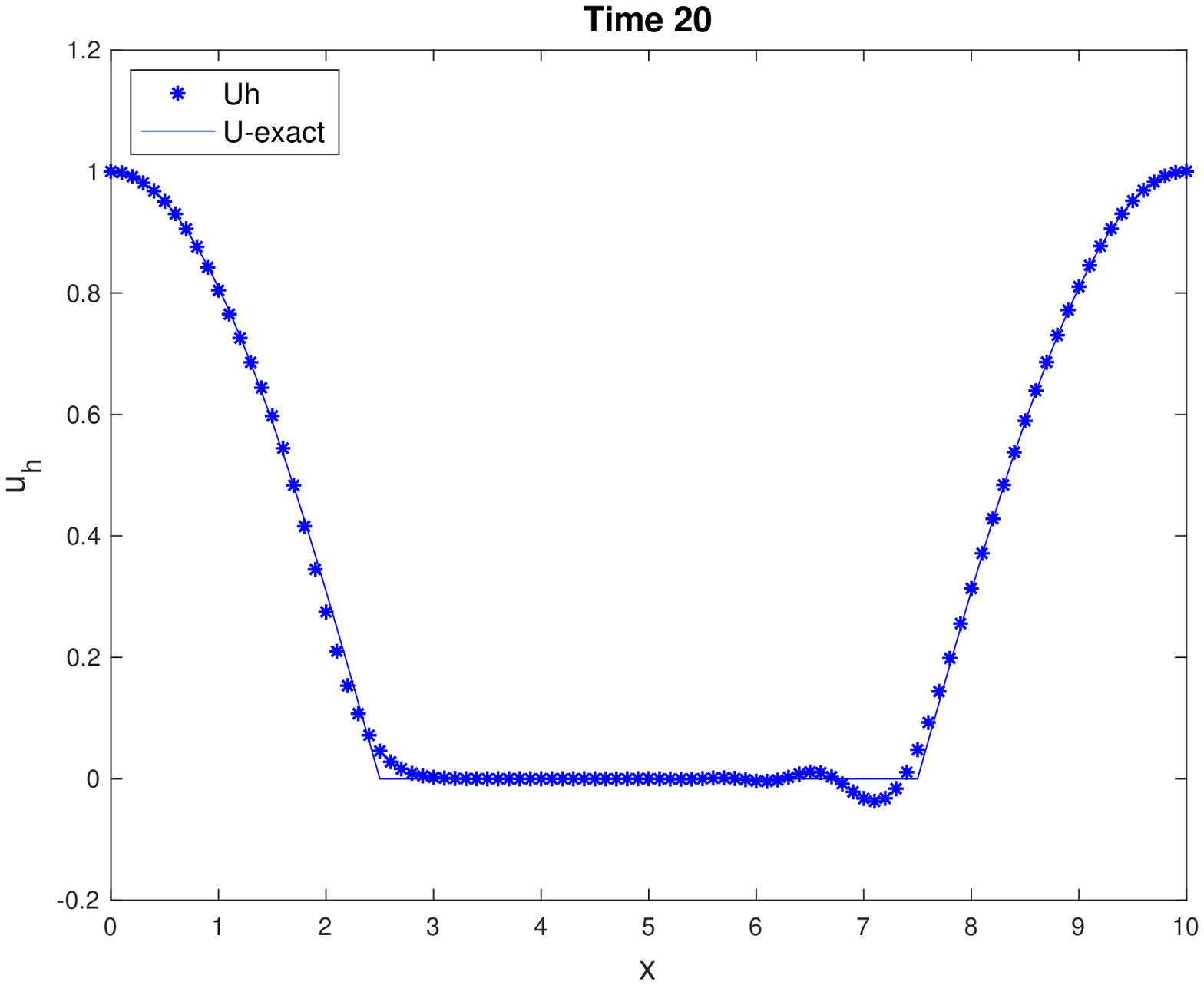}}\hspace*{8pt}\\
  \subfloat[$p=2$]{\label{Linear_Advection3} \includegraphics[width=0.48\textwidth,height=0.25\textheight]{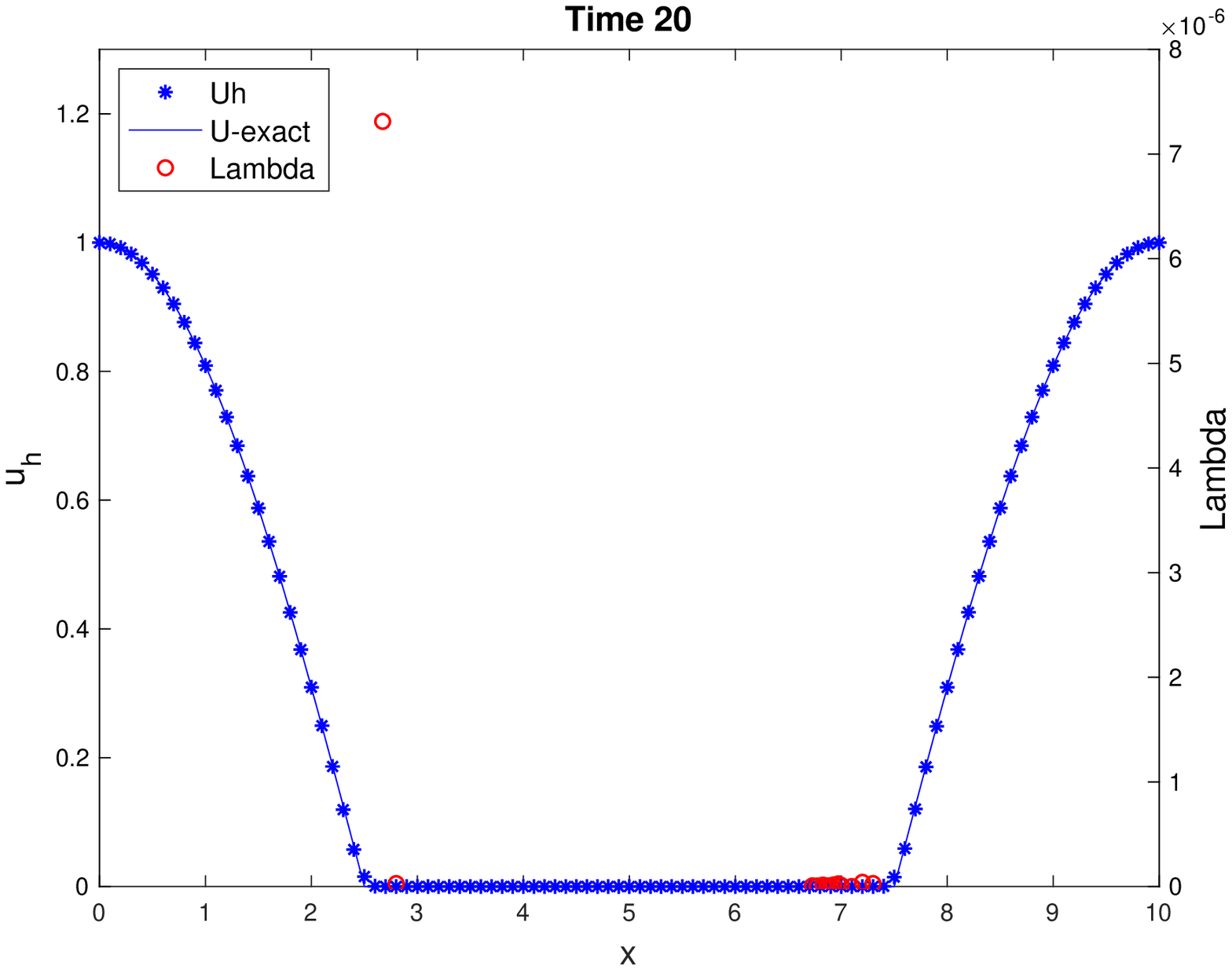}}\hspace*{5pt}
  \subfloat[$p=3$]{\label{Linear_Advection4}\includegraphics[width=0.48\textwidth,height=0.25\textheight]{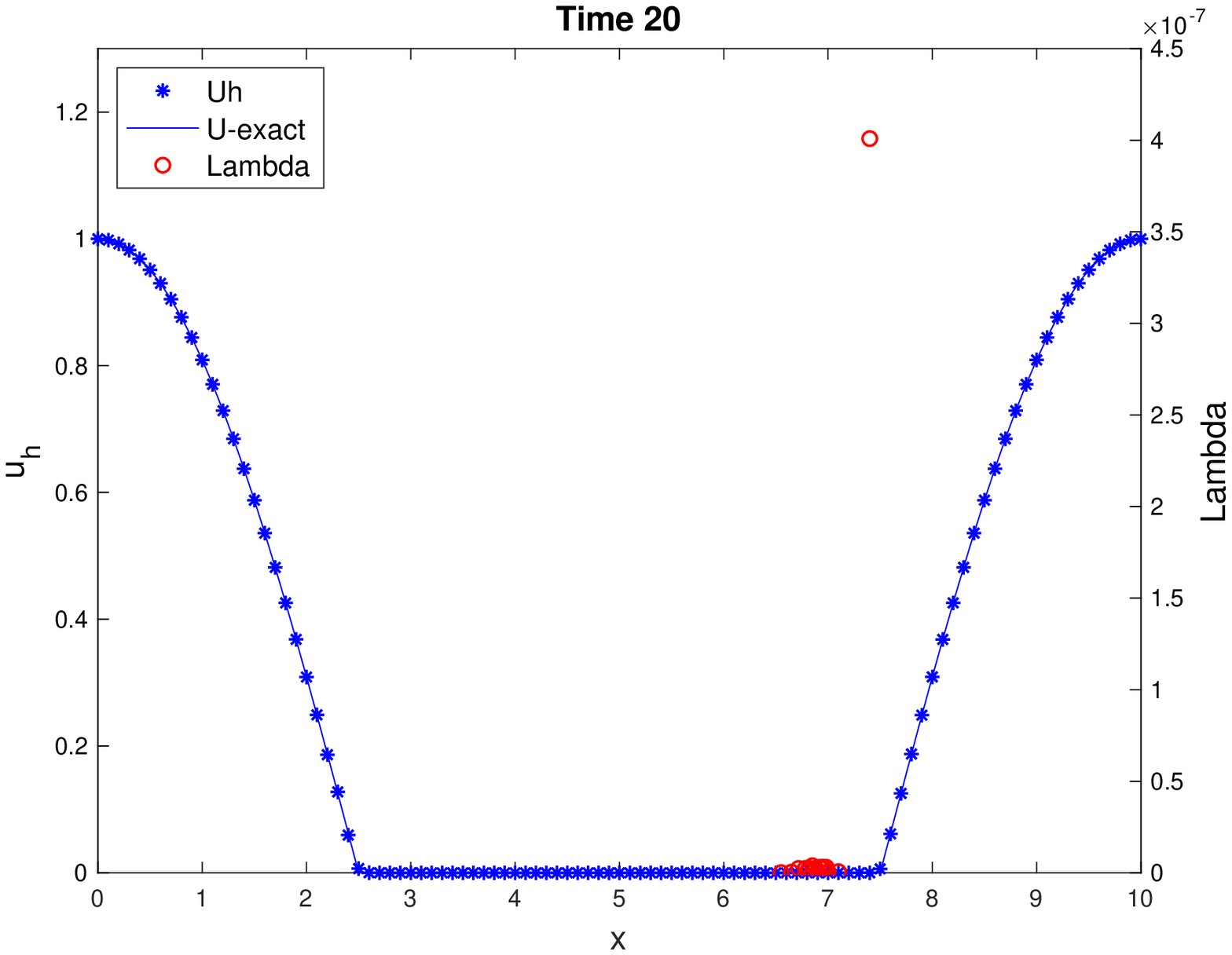}}
 \caption{Example 5.3, advection equation 1D,   (a), (c), (d) numerical solution $u_h$ with positivity preserving limiter, polynomial order, respectively, $p=1$, 2, and $3$, (b) numerical solution $u_h$ without positivity preserving limiter, polynomial order $p=1$. Computational mesh $100$ elements.  Values of the Lagrange multiplier used in the positivity preserving limiter larger than $10^{-10}$ are indicated
  in (a), (c) and (d) with a red circle.}
  \label{1DAdvection_Equation_Overview}
\end{figure}

\indent{\it Example 5.3} (1D linear advection equation). We consider \cref{linadvec_eq} with a zero right hand side in the domain $\Omega=[0,10]$ and periodic boundary conditions. The exact solution is
\begin{equation*}
u(x,t)=\max(\cos(2\pi(x-t)/10),0), \quad \text{for}\; x\in\Omega, t\in[0,T].
\end{equation*}
A constrained projection of $u(x,0)$ onto the finite element space $V_h^p$ is used as initial solution $u_h(x,0)$. The computational mesh contains 100 elements and the maximum CFL number is 1. In Figures \ref{Linear_Advection1}, \ref{Linear_Advection3},  and \ref{Linear_Advection4} the exact  and numerical solution at time $t=20$ are plotted for, respectively, polynomial orders 1, 2 and 3. At this time the wave has travelled twice through the domain and the numerical solution matches very well with the exact solution. Also, plotted is the value of the Lagrange multipliers used to impose the positivity constraint $u_{h\min}=10^{-10}$. These plots clearly show that the limiter is only active at locations where the constraint must be imposed and not in the smooth part of the solution. In Figure \ref{Linear_Advection2}, the solution for polynomial order $p=1$ without the KKT-Limiter is plotted, which clearly shows that without the limiter the solution is significantly below the $u=0$ minimum of the exact solution $u(x,t)$.

\begin{figure}[htb]
  \centering
  \subfloat[]{\label{Advection_2D_Solution}\includegraphics[width=0.49\textwidth]{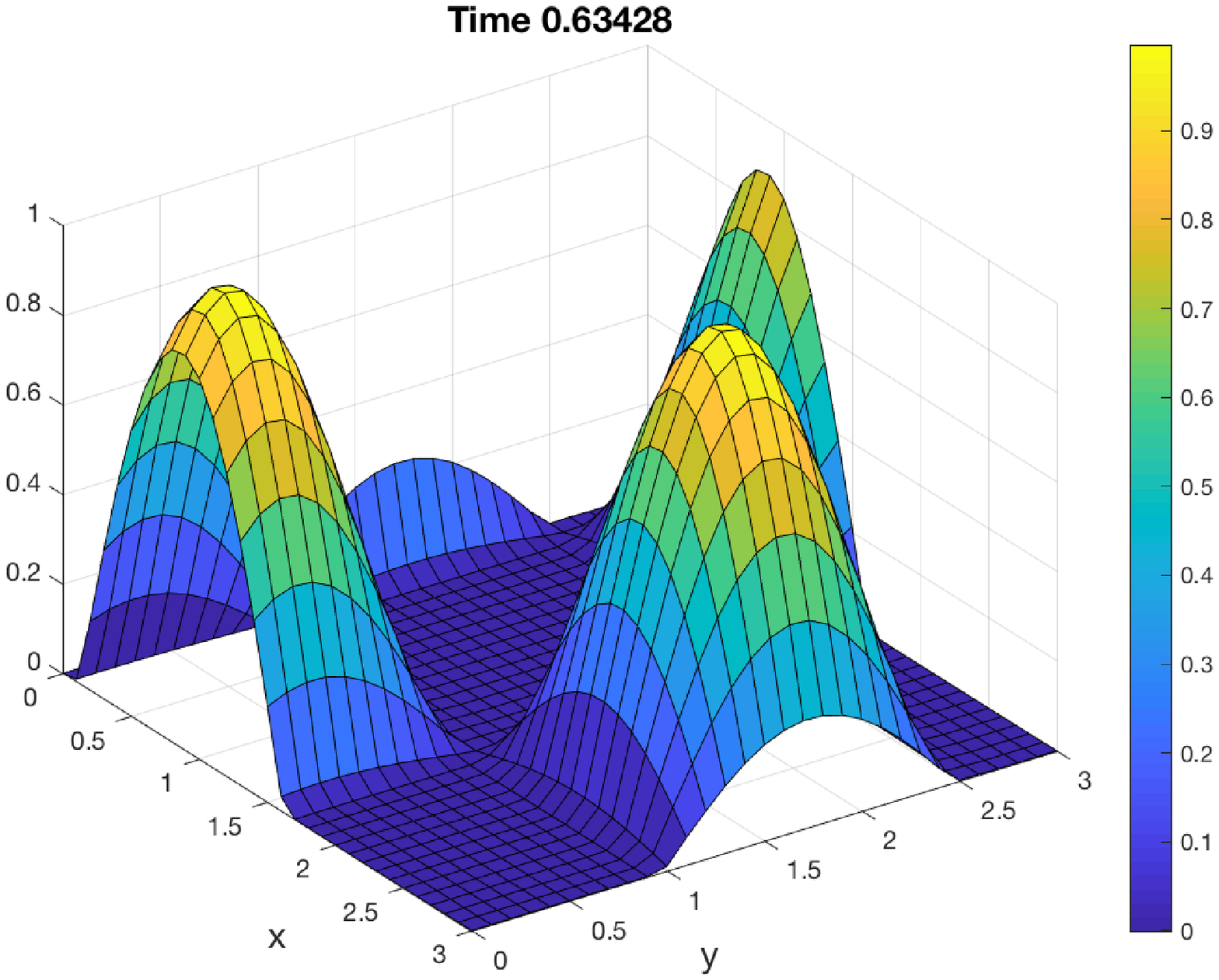}}
  \subfloat[]{\label{Advection_2D_Lagrange_Multiplier}\includegraphics[width=0.49\textwidth]{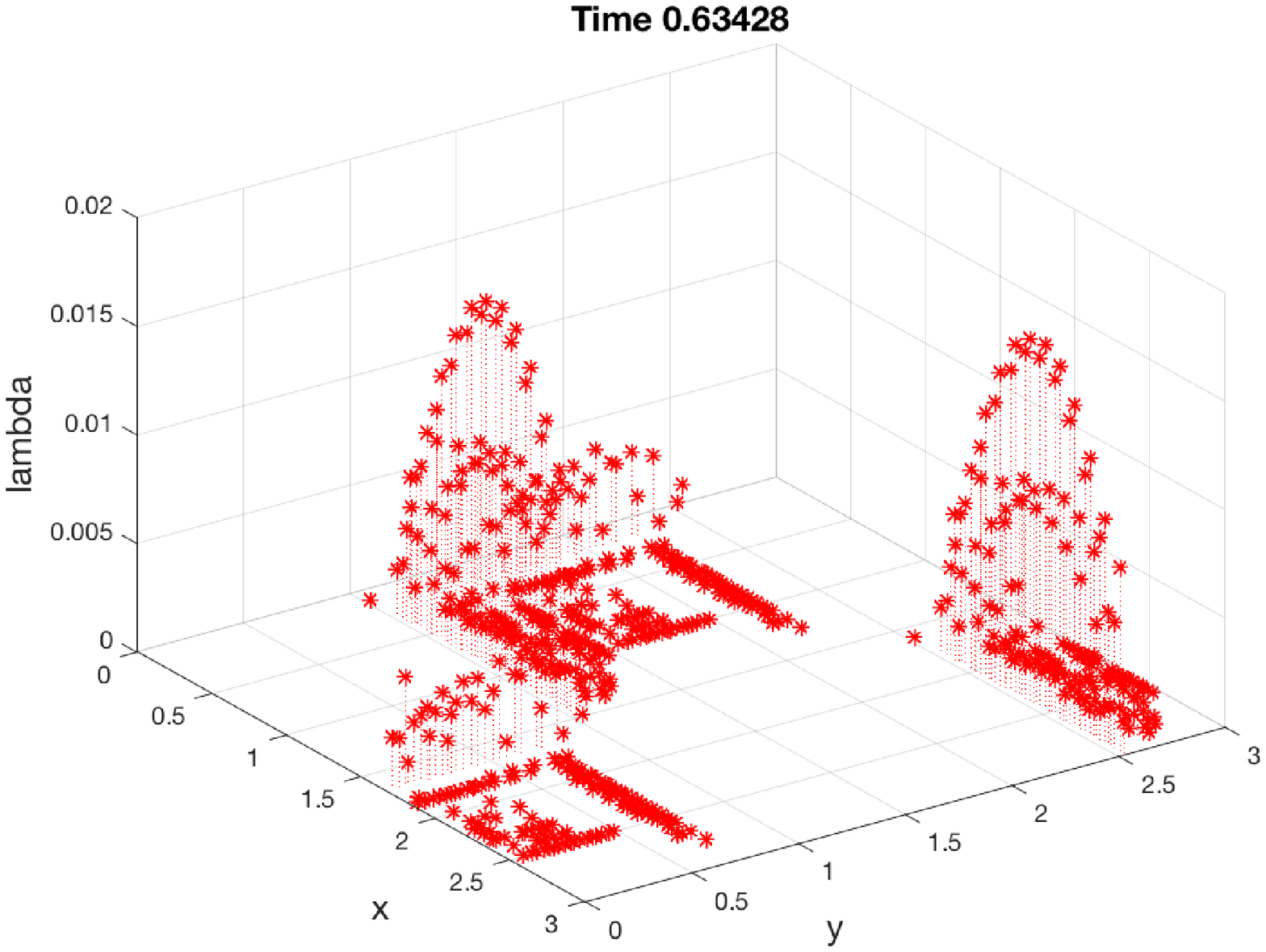}}
   \caption{Example 5.4, advection equation 2D,  (a) solution $u_h$, (b)  Lagrange multiplier. Computational mesh $30\times 30$ elements, polynomial order $p=3$. Values of the Lagrange multiplier used in the positivity preserving limiter larger than $10^{-10}$ are indicated in (b) with a red asterisk.}
  \label{2DAdvection_Equation_Overview}
\end{figure}

\indent{\it Example 5.4} (2D linear advection equation). The KKT-Limiter is also tested on a 2D linear advection equation, which is obtained by setting   $F(u)=cu$, with  $c=(-1,-2)$, and $G(u)=\nu(u)=0$ in \cref{eq:2ndconservation_law}. The domain $\Omega=[0,3]^2$ with periodic boundary conditions is used in the computations. The computational mesh contains $30\times 30$ elements. The exact solution is
\begin{equation*}
u(x,t)=\max(\cos(2\pi(x+t)/3)\cos(2\pi(y+2t)/3),0)\quad \text{for}\; x\in\Omega, t\in[0,T].
\end{equation*}
A constrained projection of $u(x,0)$ onto the finite element space $V_h^p$ is used as initial solution $u_h(x,0)$. The maximum CFL number is 1. In Figure \ref{Advection_2D_Solution} the numerical solution is shown at $t=6.3428$ and in Figure \ref{Advection_2D_Lagrange_Multiplier} the values of the Lagrange multipliers used to enforce the positivity constraint $u_{h\min}=10^{-10}$. Comparing Figures  \ref{Advection_2D_Solution} and \ref{Advection_2D_Lagrange_Multiplier} clearly shows that the KKT-Limiter is only active in those parts of the domain where the solution needs to satisfy the positivity constraint and not in the smooth part.

{ \indent{\it Example 5.5} (1D Burgers equation). 
\begin{figure}[htb]
  \centering
  \subfloat[]{\label{Burgers_1D_Solution_A}\includegraphics[width=0.49\textwidth]{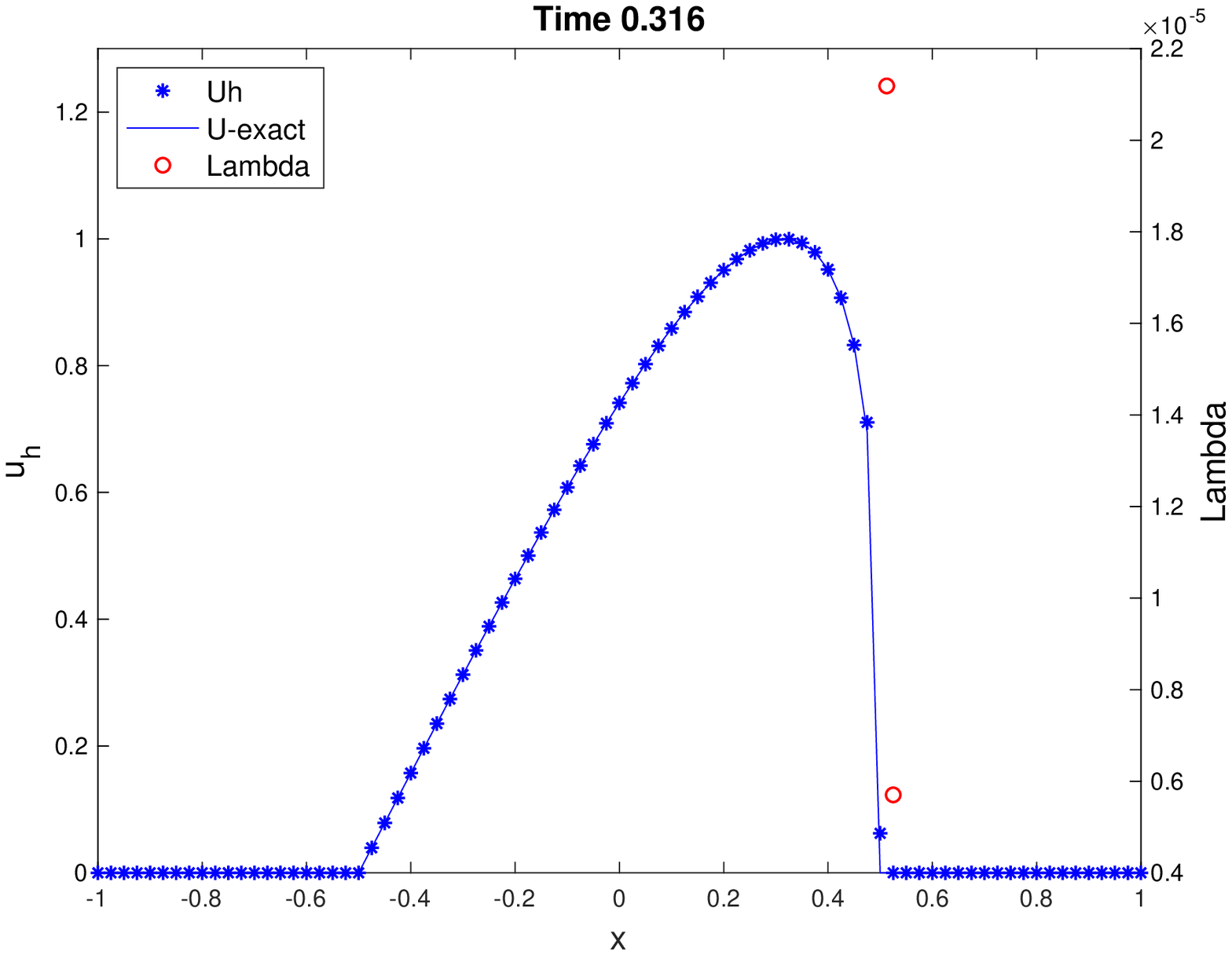}}\hspace*{3pt}
  \subfloat[]{\label{Burgers_1D_Solution_B}\includegraphics[width=0.49\textwidth]{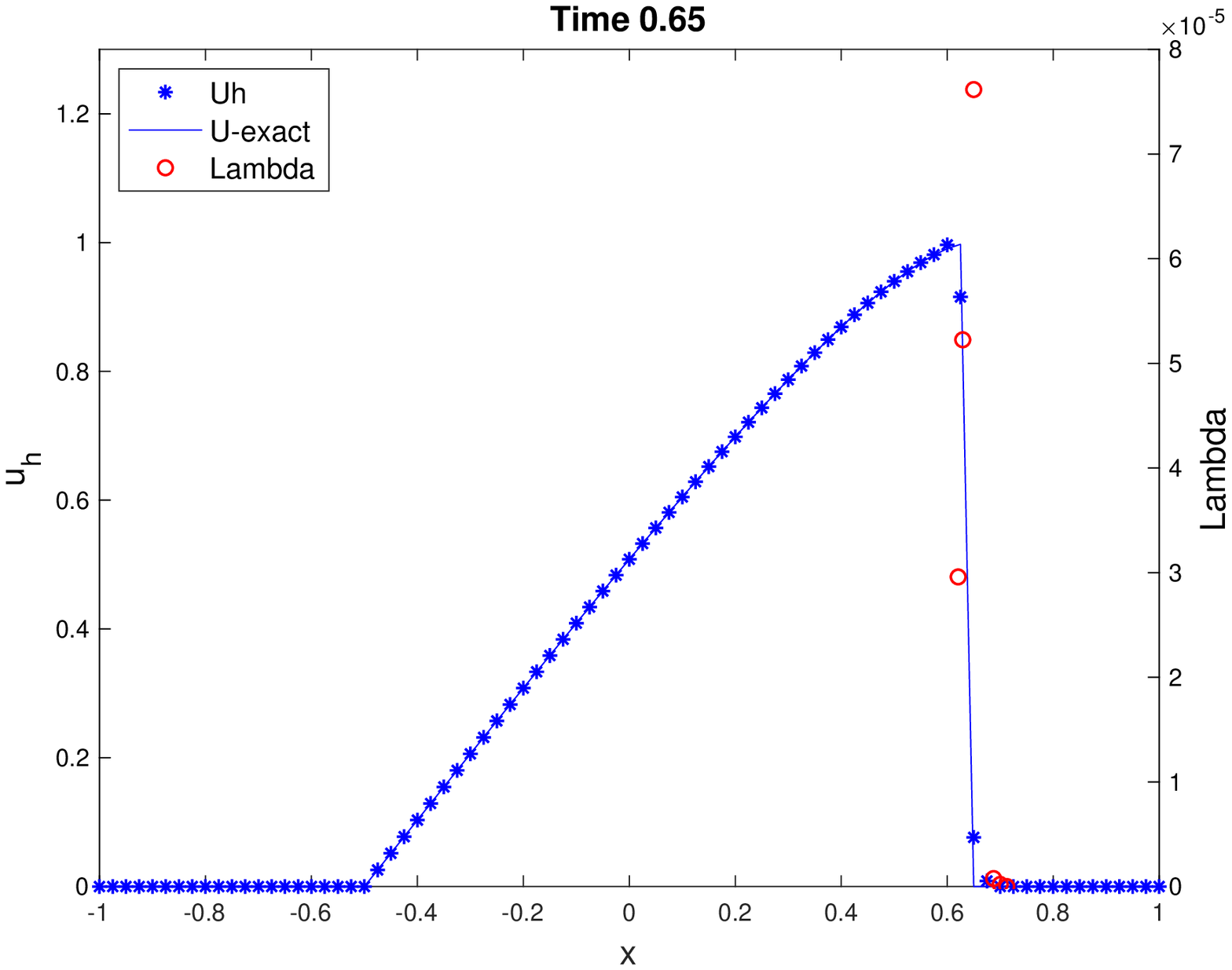}}\\
  \subfloat[]{\label{Burgers_1D_Solution_C}\includegraphics[width=0.49\textwidth]{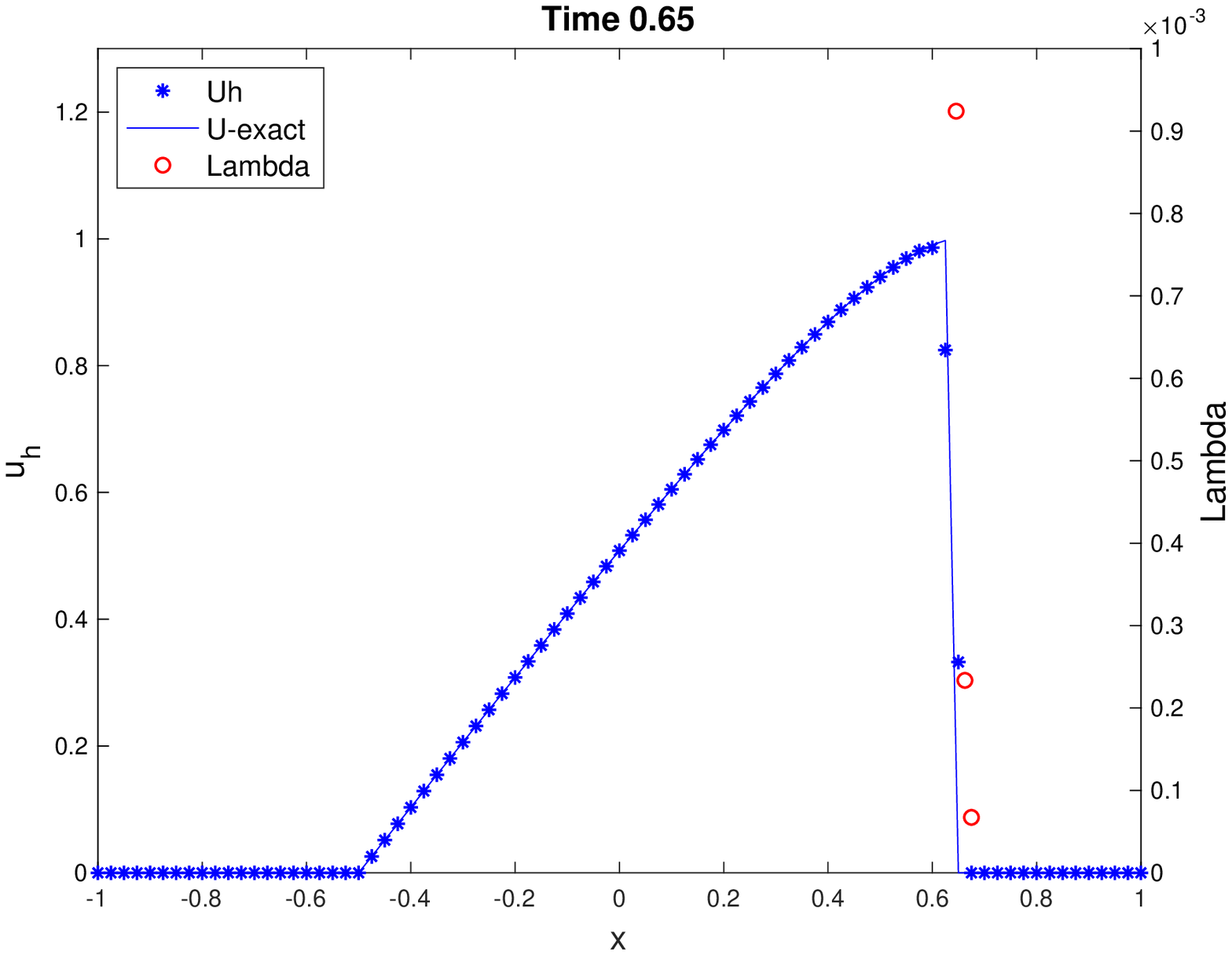}}
   \caption{Example 5.5, Burgers equation 1D,  (a)-(c) solution $u_h$ and Lagrange multiplier. The solution in (a) and (b) is computed with local conservation imposed as an explicit constraint, whereas (c) shows the solution without explicitly imposing local conservation. Computational mesh $80$ elements, polynomial order $p=3$. Values of the Lagrange multiplier used in the positivity preserving limiter larger than $10^{-10}$ are indicated in  with a red circle.}
  \label{1DBurgers_Equation_Overview}
\end{figure}
In order to test the KKT-Limiter on problems with time-dependent shocks we consider the 1D Burgers equation on a domain $\Omega=[-1,1]$ with initial condition $u_0=\max(\cos(\pi x),0)$ and periodic boundary conditions. The polynomial order is $p=3$. As lower  and upper bounds in the positivity preserving limiter we use, respectively $u_{h\min}=10^{-10}$ and $u_{h\max}=1$, and no monotonicity constraint is imposed. The initially smooth part of the solution develops into a shock. The onset of the shock is shown in Figure \ref{Burgers_1D_Solution_A} and the later stages of the shock at $t=0.65$ in Figure \ref{Burgers_1D_Solution_B}. Figure \ref{Burgers_1D_Solution_C} shows the solution when the conservation constraint \cref{cons_constr} is not explicitly enforced. The difference in the shock solution for the discretizations with and without the explicitly imposed conservation constraint is very small. The main reason for this is that the KKT-Limiter is only active in regions where the constraints must be imposed and does not affect the discretization at other places in the domain. This can be seen from the values of the Lagrange multipliers that are used to impose the positivity constraints, which are indicated with red circles, and are only non-zero in the vicinity of the shock and at locations where the solution has a discontinuous derivative. The KKT-Limiter to ensure the positivity constraints therefore has a very small effect on the conservation properties of the DG discretization as can be seen by comparing Figures \ref{Burgers_1D_Solution_B} and \ref{Burgers_1D_Solution_C}.}

\indent{\it Example 5.6} (Allen-Cahn equation). The Allen-Cahn equation is a reaction-diffusion equation that describes phase transition. The Allen-Cahn equation is obtained by setting  $G(u)=u^3-u$,  $\nu(u)=\bar{\nu}$, and $F(u)=0$ in \cref{eq:2ndconservation_law}. The solution of the Allen-Cahn equation should stay within the range $[0,1]$. Hence, we apply both the positivity  and maximum preserving limiters, respectively, \cref{positivity_limiter}-\cref{maximum_limiter} with bounds $u_{h{\rm min}}=10^{-14}$ and $u_{h{\rm max}}=1-10^{-10}$. A constrained projection of $u(x,0)$ onto the finite element space $V_h^p$ is used as initial solution $u_h(x,0)$.

\begin{figure}[tbhp]
  \centering
  \subfloat[$u_h$ - with limiter]{\label{Allen_Cahn1}\includegraphics[width=0.49\textwidth,height=0.22\textheight]{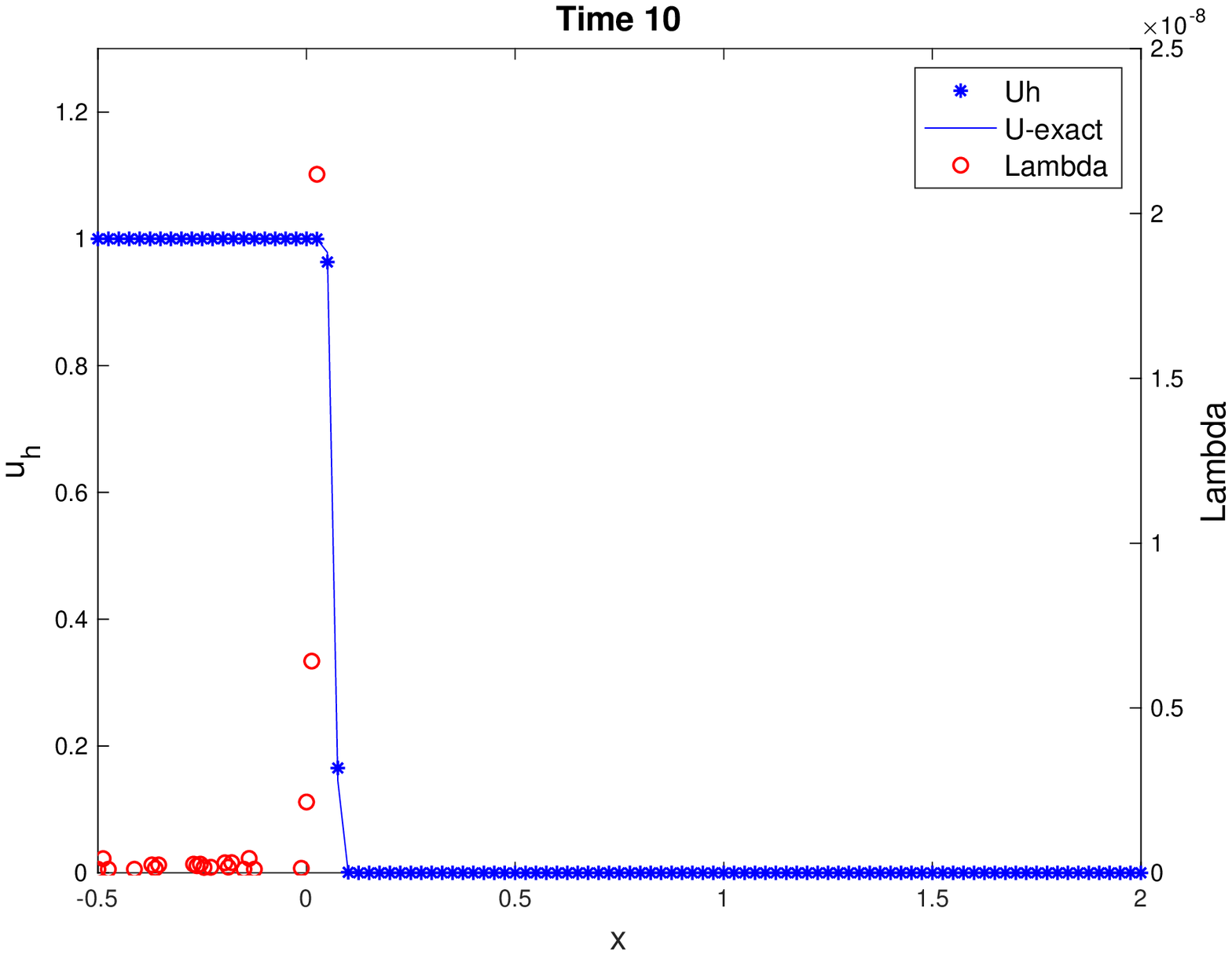}}\hspace*{5pt}
  \subfloat[$Q_h$  - with limiter]{\label{Allen_Cahn2}\includegraphics[width=0.49\textwidth,height=0.22\textheight]{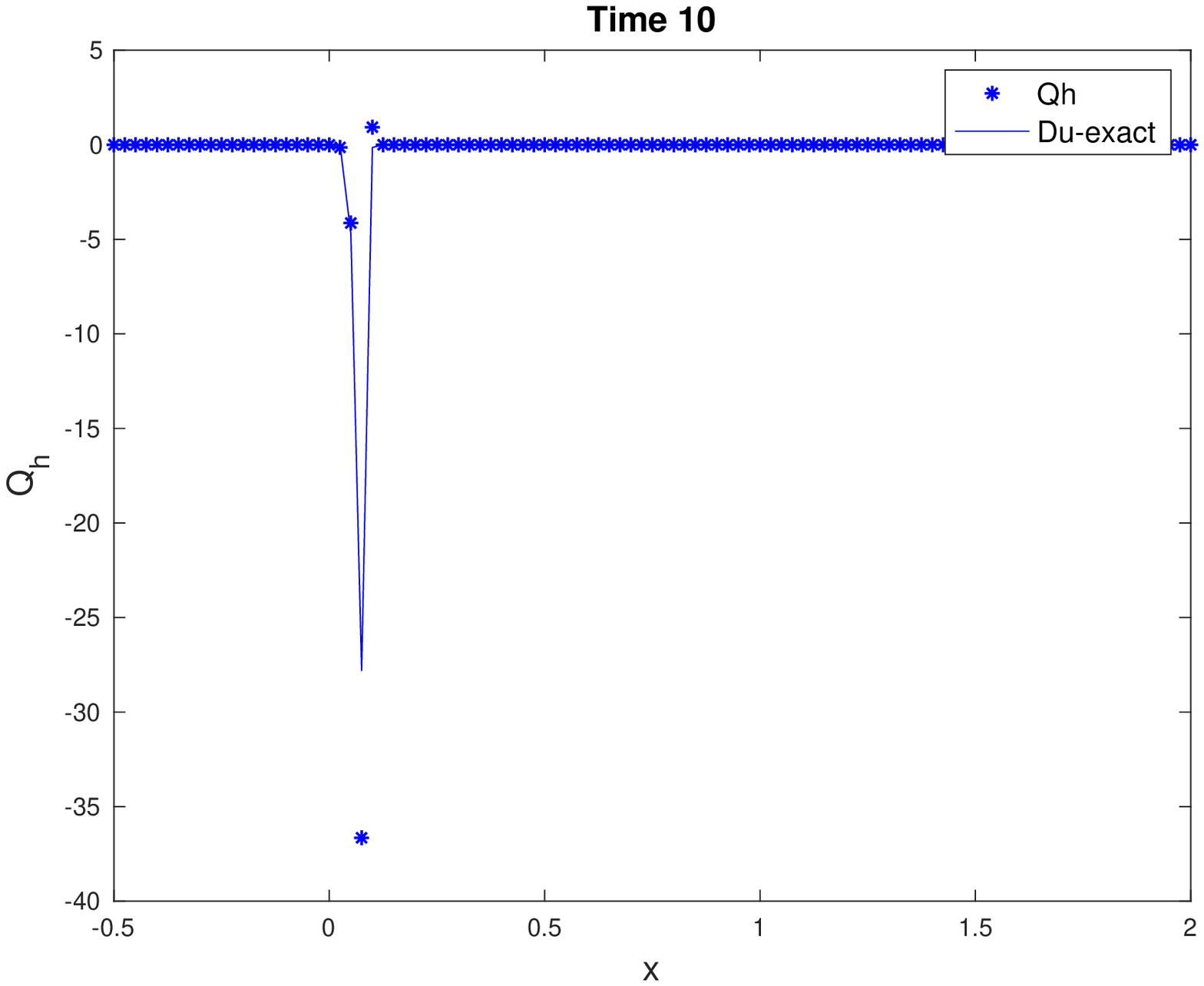}}
 \caption{Allen-Cahn equation 1D, Example 5.6a, (a) numerical solution $u_h$ and exact solution $u$, (b) derivative of numerical solution  $Q_h$ and exact derivative $Du$.  Computational mesh $100$ elements, polynomial order $p=3$.  Values of the Lagrange multiplier used in the positivity and maximum preserving limiters larger than $10^{-10}$ are indicated in (a) with a red circle.}
  \label{Allen_Cahn_Overview1}
\end{figure}

\indent{\it Example 5.6a} (Allen-Cahn equation 1D).  As test case we use the traveling wave solution
\begin{equation*}
u(x,t)= \frac{1}{2}\left(1-\tanh\left(\frac{x-st}{2\sqrt{2\bar{\nu}}}\right)\right),
\end{equation*}
with wave velocity $s=3\sqrt{\bar{\nu}/2}$. The computational domain is $\Omega = [-\frac{1}{2},2]$. If the mesh resolution is sufficiently dense such that the jump in the traveling wave solution is well resolved, then no limiter is required. For  small values of the viscosity the solution will, however, violate the positivity constraints, except on very fine meshes. In Figures  \ref{Allen_Cahn1}  and \ref{Allen_Cahn2}, respectively, the numerical solution $u_h$ and its derivative $Q_h$ and the exact solutions are shown for the viscosity $\bar{\nu}=10^{-5}$ on a mesh with 100 elements and polynomial order 3 for the basis functions. The values of the Lagrange multiplier used to impose the positivity constraints are also shown in Figure \ref{Allen_Cahn1}. The solution has a very thin and steep transition region, but the wave speed is still correctly computed by the LDG scheme and the KKT limiter ensures that both the positivity  and  maximum constraint are satisfied.

\begin{figure}[tbhp]
  \centering
  \subfloat[]{\label{Allen_Cahn3}\includegraphics[width=0.49\textwidth,height=0.25\textheight]{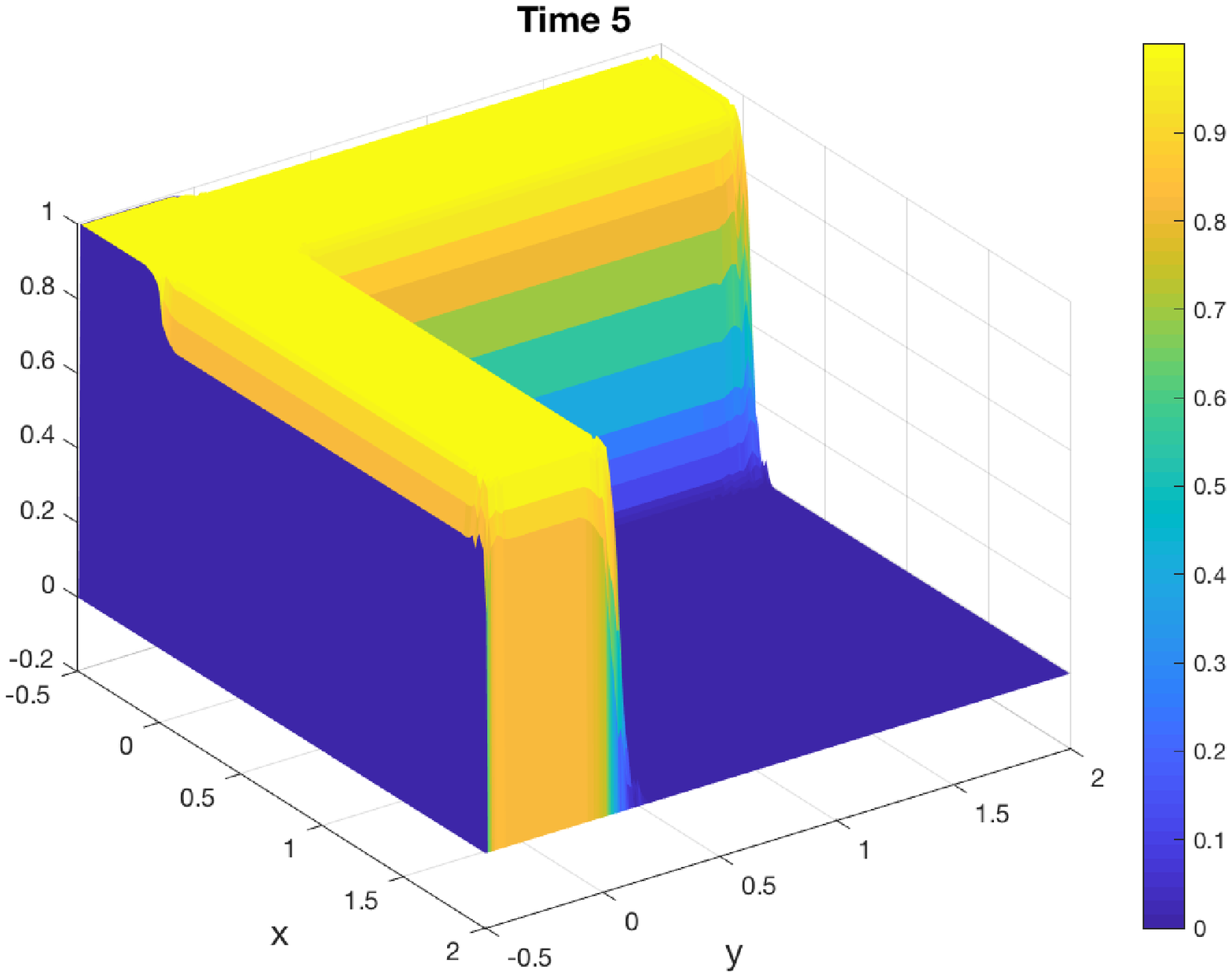}}
  \subfloat[]{\label{Allen_Cahn4}\includegraphics[width=0.49\textwidth,height=0.25\textheight]{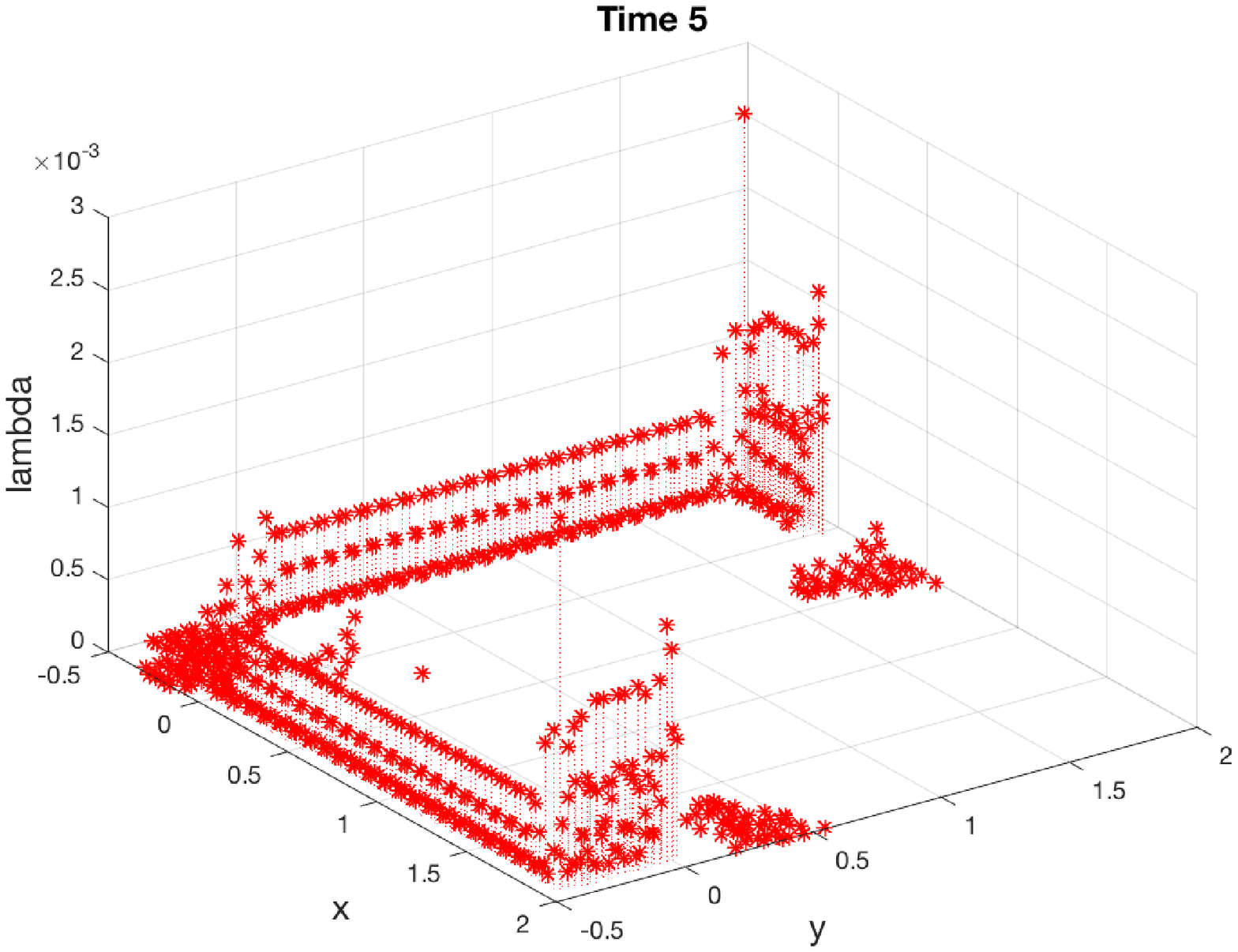}}
 \caption{Allen-Cahn equation 2D, Example 5.6b, (a) numerical solution $u_h$ (b)  Lagrange multiplier. Computational mesh $30\times 30$ elements, polynomial order $p=3$.  Values of the Lagrange multiplier used in the positivity and maximum preserving limiters larger than $10^{-10}$ are indicated in (b) with a red asterisk.}
  \label{Allen_Cahn_Overview2}
\end{figure}

\indent{\it Example 5.6b} (Allen-Cahn equation 2D). For the 2D test case the computational domain is $\Omega = [-\frac{1}{2},2]^2$ and the computational mesh contains $30\times 30$ elements. The viscosity coefficient is selected as $\bar{\nu}=10^{-4}$. As test case we use the initial solution
\begin{equation*}
u(x,0)= \frac{1}{4}\left(1-\tanh\left(\frac{x}{2\sqrt{2\bar{\nu}}}\right)\right)\left(1-\tanh\left(\frac{y}{2\sqrt{2\bar{\nu}}}\right)\right),
\end{equation*}
which values are also used as boundary condition for $t>0$.
At this mesh resolution a positivity preserving limiter is necessary. The numerical solution shown in Figure \ref{Allen_Cahn3} has steep gradients and the positivity preserving limiter  ensures that the bounds are satisfied. The locations where the limiter is active can be seen in Figure \ref{Allen_Cahn4}, which shows the values and locations of the Lagrange multipliers used to impose the bounds in the DG discretization.


\indent{\it Example 5.7} (Barenblatt equation). The Barenblatt equation, which models a porous medium, is obtained by setting $\nu(u)=mu^{m-1}$, $m>1$, and $F(u)=0$, $G(u)=0$ in \cref{eq:2ndconservation_law}.
The exact solution is
\begin{equation*}
u(t,x)=t^\alpha\left(\left(C-\frac{\beta(m-1)}{2m}\frac{\vert x\vert^2}{t^{2\beta}}\right)_+\right)^{\frac{1}{m-1}},
\end{equation*}
with $\alpha=\frac{n}{n(m-1)+2}$, $\beta=\frac{\alpha}{n}$,  $n={\rm dim}(\Omega)$, $(x)_+=\max(x,0)$ and $C>0$. We selected $C=1$ and $m=8$. The solution should be positive or zero for $t>0$. The initial solution for the computations is the constrained projection of $u(x,1)$ onto the finite element space $V_h^p$. In the computations Dirichlet boundary conditions are imposed, where the solution for $t>0$ is fixed at the same level as  the initial solution.

\begin{figure}[htb]
  \centering
 \includegraphics[width=0.45\textwidth]{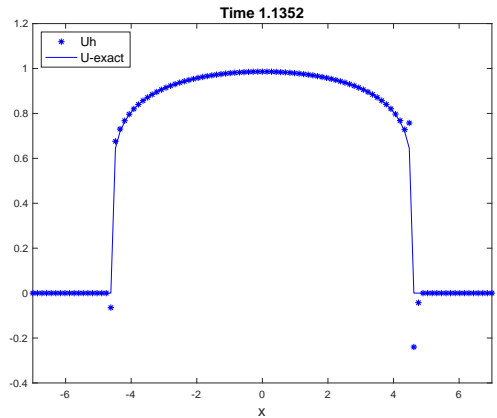}
   \caption{Barenblatt equation 1D, Example 5.7a, numerical solution $u_h$  without limiter and exact solution $u$. Computational mesh $100$ elements, polynomial order $p=3$. }
    \label{Barenblatt1D_nolimiter}
\end{figure}

\indent{\it Example 5.7a} (1D Barenblatt equation).
We first consider the 1D Barenblatt equation on the domain $\Omega=[-7,7]$ using a computational mesh of 100 elements. In Figure \ref{Barenblatt1D_nolimiter} the numerical solution without the use of a limiter is shown. It is clear that near the boundary of $u(t,x)>0$, where the derivative of $u$ becomes unbounded, significant negative values of $u_h$ are obtained. These cause severe numerical problems and do not allow the continuation of the computations.

\begin{figure}[tbhp]
  \centering
  \subfloat[]{\label{Barenblatt1}\includegraphics[width=0.49\textwidth]{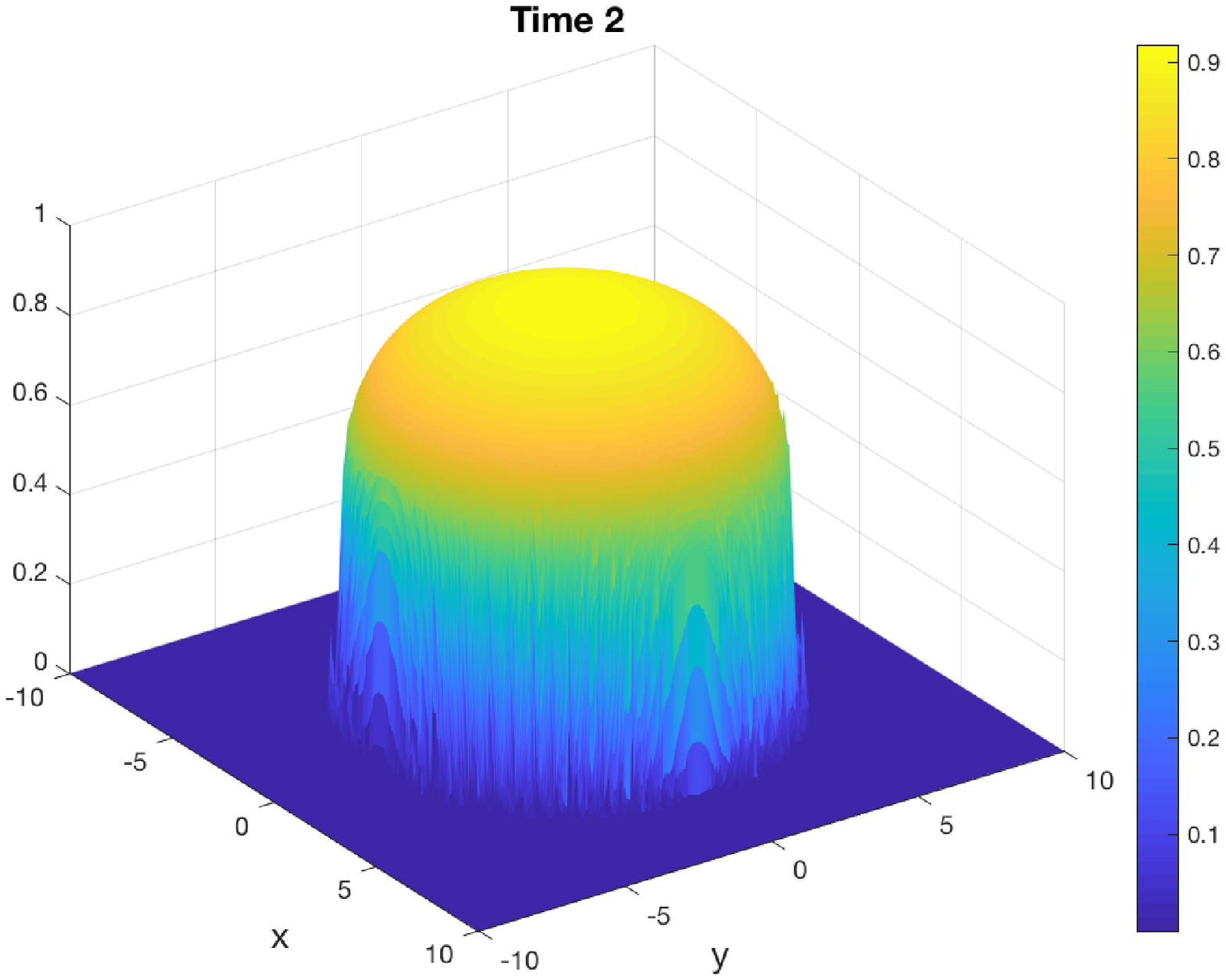}}
  \subfloat[]{\label{Barenblatt2}\includegraphics[width=0.45\textwidth]{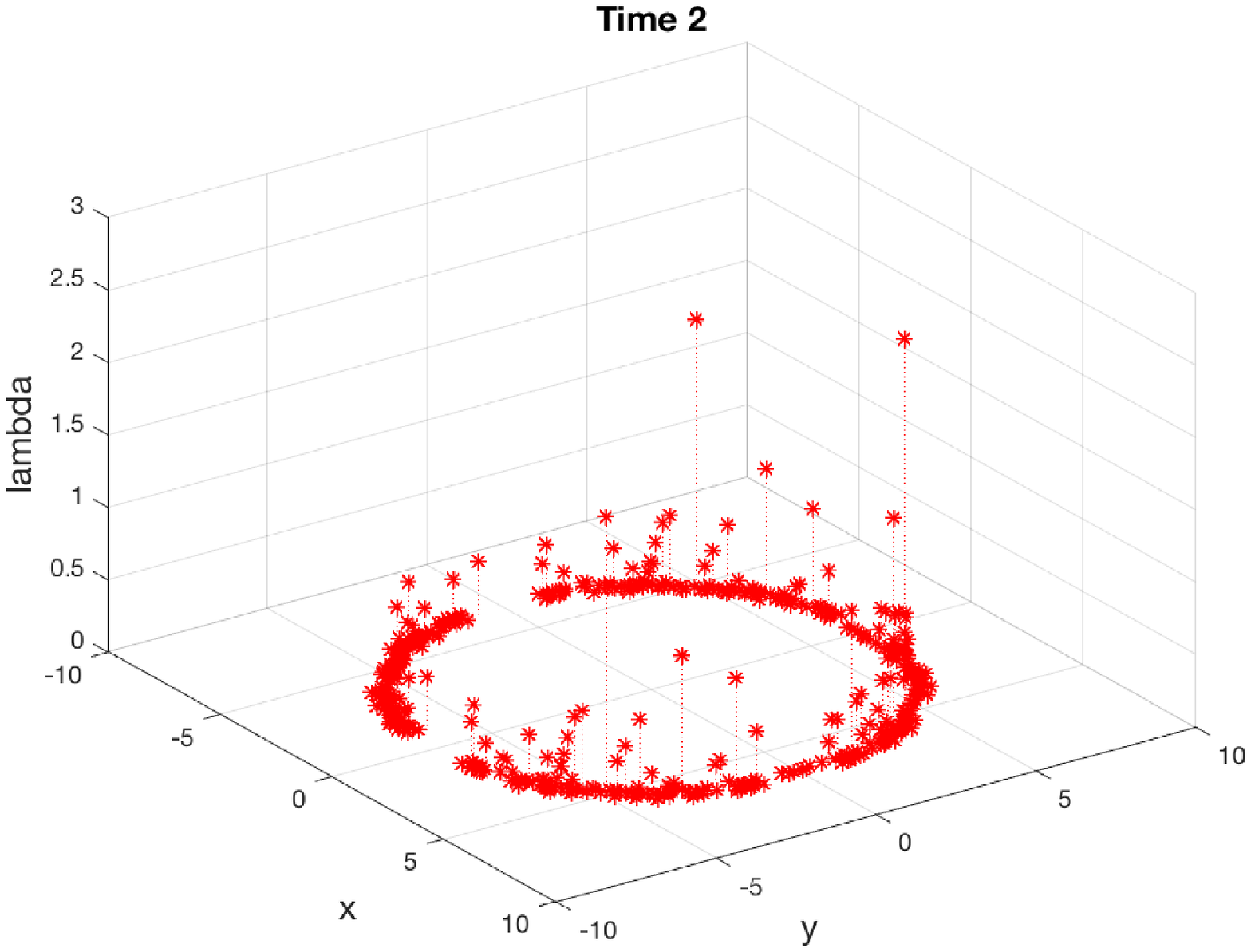}}\\[5pt]
  \subfloat[]{\label{Barenblatt3}\includegraphics[width=0.45\textwidth]{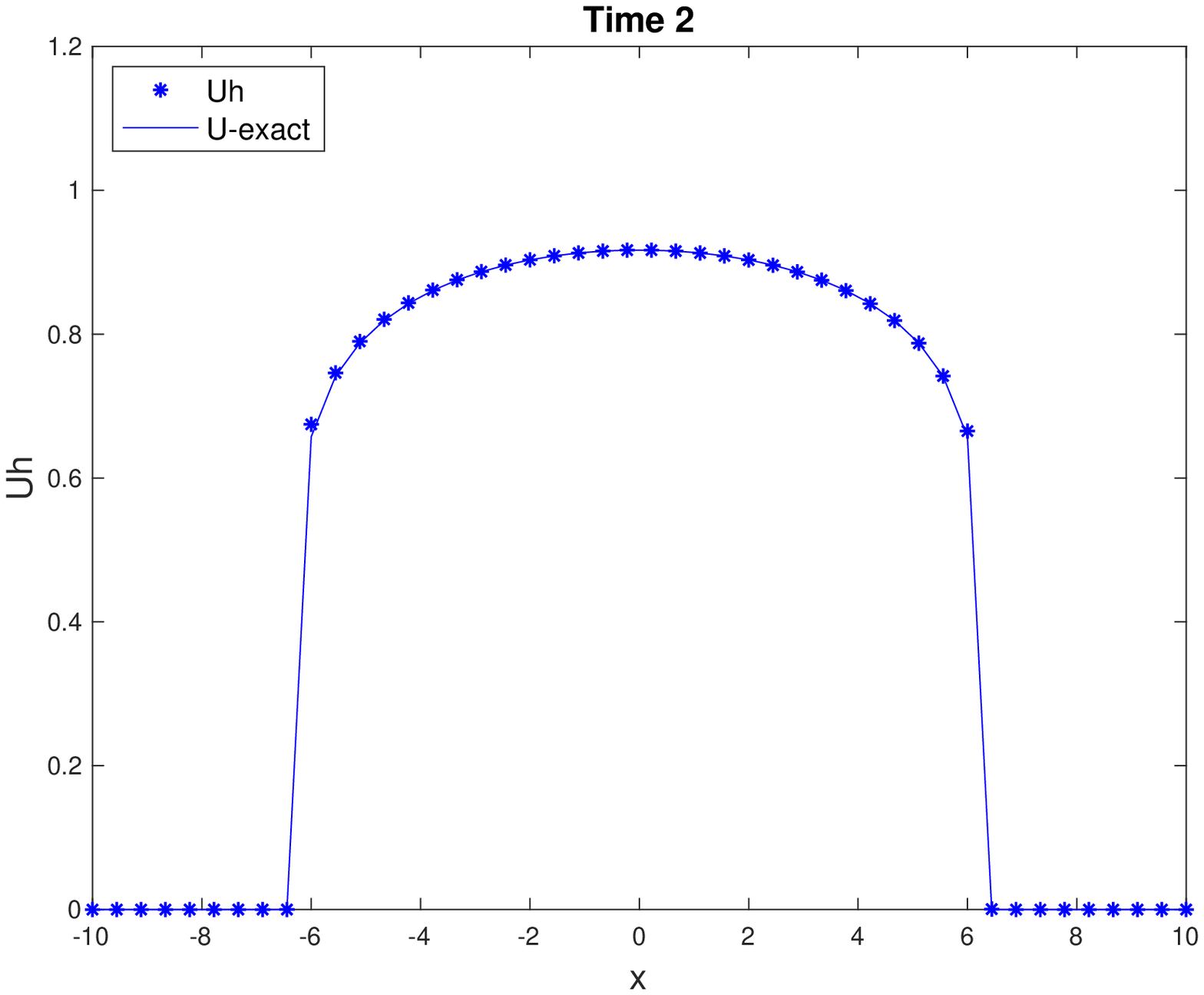}}
   \caption{Barenblatt equation 2D, Example 5.7b, (a) solution $u_h$, (b)  Lagrange multiplier, (c) numerical solution $u_h$ and exact solution $u$ in cross-section at $y=0$. Computational mesh $50\times 50$ elements, polynomial order $p=3$. Values of the Lagrange multiplier used in the positivity preserving limiter larger than $10^{-10}$ are indicated in (b) with a red asterisk.}
  \label{Barenblatt_Overview}
\end{figure}

\indent{\it Example 5.7b} (2D Barenblatt equation).
In Figures \ref{Barenblatt1} and \ref{Barenblatt2}, respectively, the numerical solution $u_h$ of the 2D Barenblatt equation and the values of the Lagrange multiplier are shown at time $t=2$ on a mesh of $50\times 50$ elements. In these computations the KKT Limiter was used, which successfully prevents the numerical solution $u_h$ from becoming negative, which is shown in Figure \ref{Barenblatt3}. The imposed constraint is $u_{h{\rm min}}=10^{-10}$. Figure \ref{Barenblatt3} also shows an excellent agreement between the exact solution $u$ and the numerical solution $u_h$.

\indent{\it Example 5.8} (1D Buckley-Leverett equation). The Buckley-Leverett equation models two phase flow in a porous medium. We consider two cases, respectively, with and without gravity. Since the solution has to be strictly inside the range $[0,1]$ we use both the positivity and maximum preserving limiter, with bounds $u_{h{\rm min}}=10^{-10}$ and $u_{h{\rm max}}=1-10^{-10}$, respectively. The computational domain is $\Omega =[0,1]$. A Dirichlet boundary condition at $x=0$, based on the initial solution, and an outflow boundary condition at $x=1$ are imposed. The viscosity coefficient is $\bar{\nu}=0.01$. Since we do not have an exact solution to compare with we compute the numerical solution on two meshes, viz. with 100 and 200 elements.
The two test cases given by Examples 5.8a and 5.8b are also considered in \cite{kurganov2000new}.

\begin{figure}[tbhp]
  \centering
  \subfloat[$u_h$ - no gravity]{\label{Buckley_Leverett1}\includegraphics[width=0.49\textwidth,height=0.25\textheight]{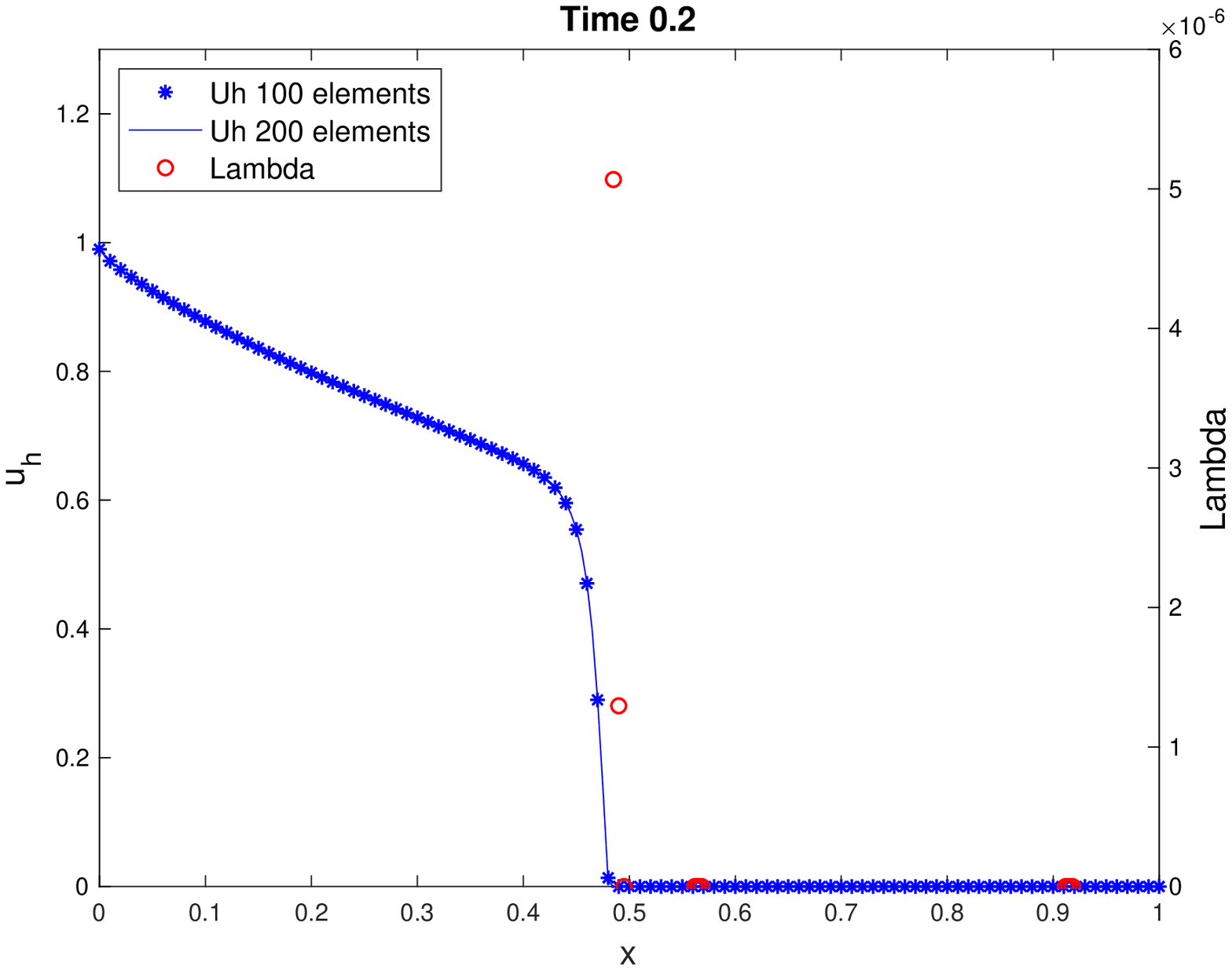}}\hspace*{5pt}
  \subfloat[$Q_h$  - no gravity]{\label{Buckley_Leverett2}\includegraphics[width=0.49\textwidth,height=0.25\textheight]{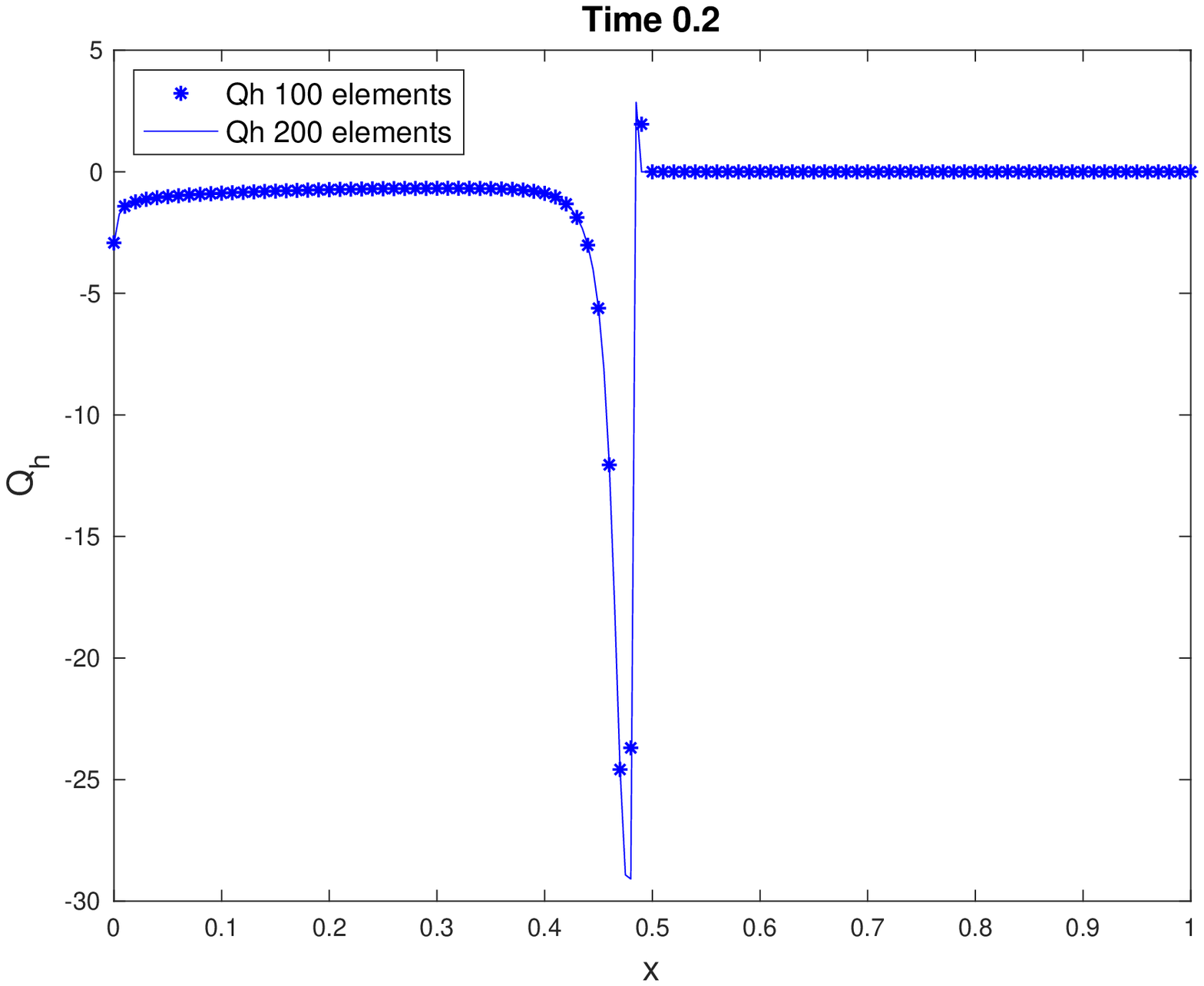}}\\[5pt]
  \subfloat[$u_h$ - gravity]{\label{Buckley_Leverett3}\includegraphics[width=0.49\textwidth,height=0.25\textheight]{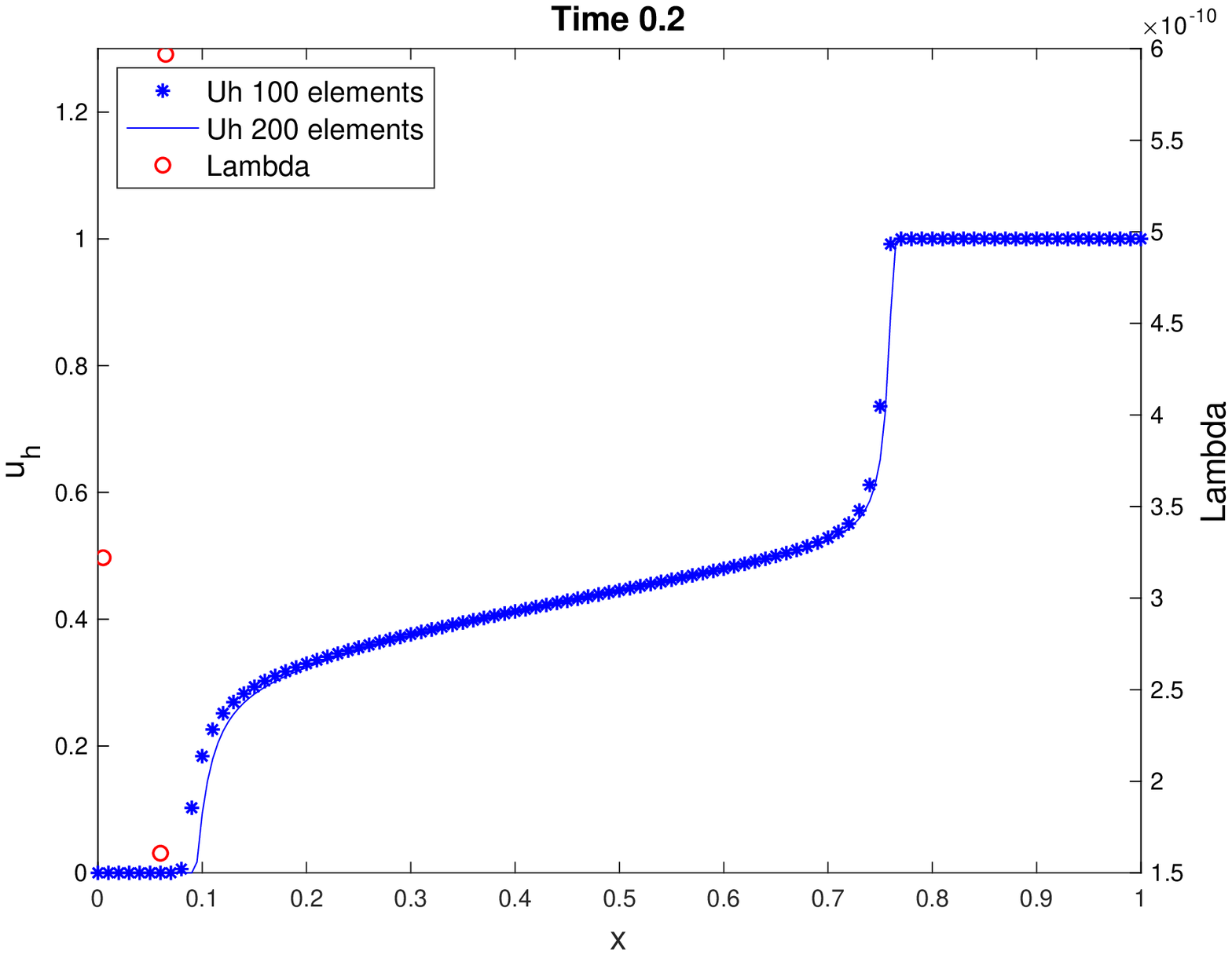}}\hspace*{5pt}
  \subfloat[$Q_h$ - gravity]{\label{Buckley_Leverett4}\includegraphics[width=0.49\textwidth,height=0.25\textheight]{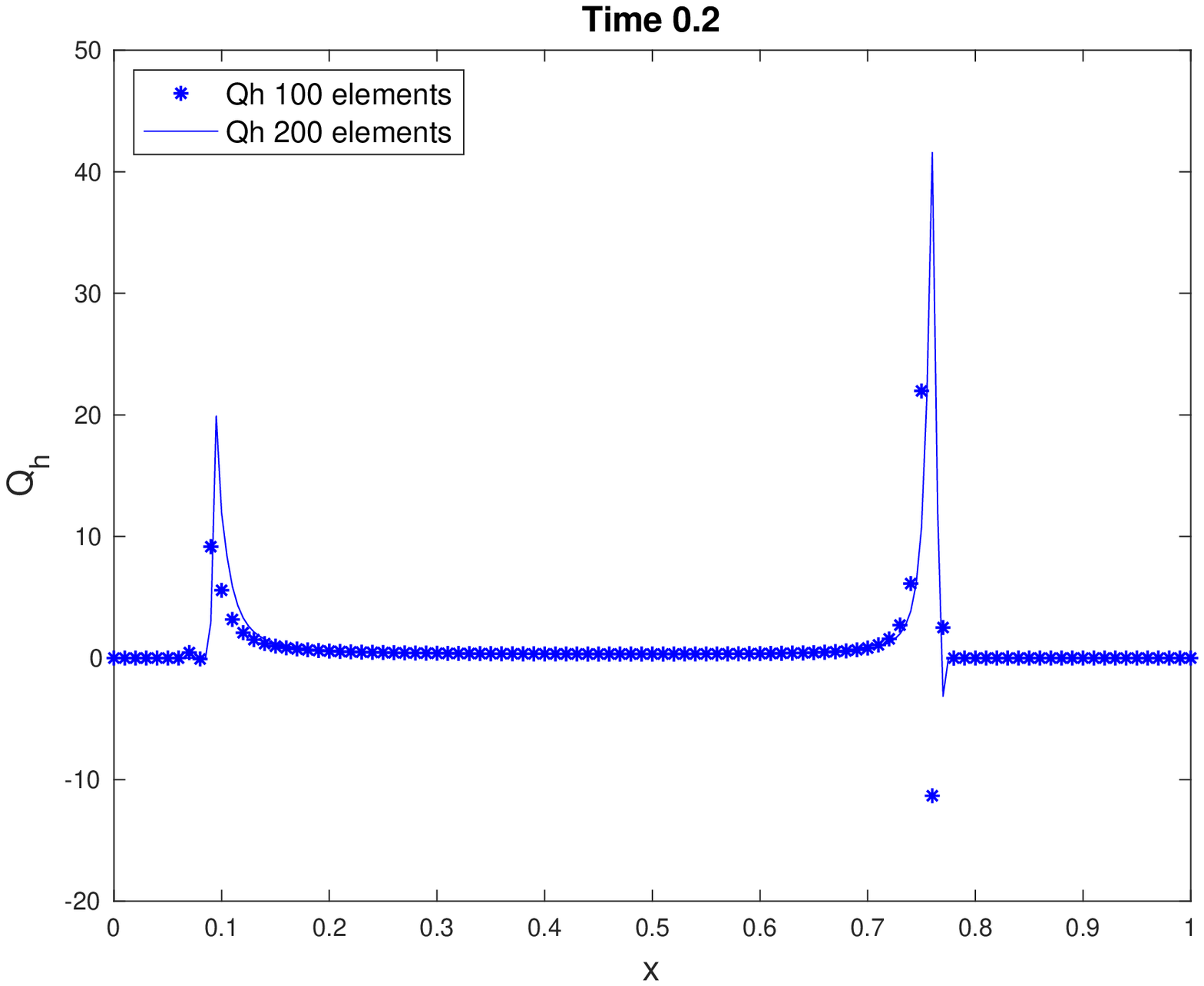}}
  \caption{Example 5.8a, Buckley-Leverett equation without gravity 1D,  (a) numerical solution $u_h$, (b) numerical solution derivative $Q_h$; Example 5.8b, Buckley-Leverett equation with gravity 1D,   (c) numerical solution $u_h$, (d) numerical solution derivative $Q_h$. Computational meshes $100$ and $200$ elements, polynomial order $p=3$. Values of the Lagrange multiplier used in the positivity and maximum preserving limiters larger than $10^{-10}$ are indicated with a red circle in (a) and (c).}
  \label{Buckley_Leverett_Overview}
\end{figure}

\indent{\it Example 5.8a} (1D Buckley-Leverett equation without gravity). The 1D Buckley-Leverett equation without gravity is obtained by setting $G(u)=0$, and $\nu(u)$ and $F(u)=f(u)$, respectively, as
\begin{equation*}
\nu(u)=
\begin{cases}
4\bar{\nu} u(1-u),\hspace*{10pt}&\text{if}\; 0\leq u\leq 1,\\
0,&\text{otherwise}.
\end{cases}
\end{equation*}
\begin{equation}
f(u)=
\begin{cases}
0,&\text{if}\; u<0,\\
\frac{u^2}{u^2+(1-u)^2},\hspace*{6pt}&\text{if}\; 0\leq u\leq 1,\\
1,&\text{if}\; u>1.
\end{cases}\label{BLflux_nogravity}
\end{equation}
The initial condition is
\begin{equation*}
u(x,0)=\begin{cases}
0.99-3x\quad& 0\leq x\leq 0.33,\\
0&\frac{1}{3} <x\leq 1.
\end{cases}
\end{equation*}
The numerical solution $u_h$ and its derivative $Q_h$ are shown in, respectively, Figures \ref{Buckley_Leverett1} and \ref{Buckley_Leverett2}.  Also, the values of the Lagrange multiplier used to enforce the constraints is shown in Figure \ref{Buckley_Leverett1}. The limiter is only active in the thin layer between the phases and is crucial to obtain sensible physical solutions.  The results of 100 and 200 elements match well.

\indent{\it Example 5.8b} (1D Buckley-Leverett equation with gravity). A much more difficult test case is provided by the Buckley-Leverett equation with gravity, which is obtained by modifying the flux $F(u)$ as
\begin{equation*}
F(u)=\begin{cases}
f(u)(1-5(1-u)^2),\quad&  u\leq 1,\\
1& u>1,
\end{cases}
\end{equation*}
with $f(u)$ given by \cref{BLflux_nogravity}. The initial solution is
\begin{equation*}
u(x,0)=\begin{cases}
0\quad& 0\leq x\leq a,\\
\frac{1}{mh}(x-a)& a< x \leq 1-\frac{1}{\sqrt{2}},\\
1&1-\frac{1}{\sqrt{2}} <x\leq 1,
\end{cases}
\end{equation*}
with $a=1-\frac{1}{\sqrt{2}}-mh$, $h$ the mesh size and $m=3$. The linear transition for $x$ in the range $[a,1-\frac{1}{\sqrt{2}}]$ is used to remove the infinite value in the derivative, which would otherwise result in unbounded values of $Q_h$ at $t=0$.
The Buckley-Leverett equations with gravity result a strongly nonlinear problem where the equations change type and is a severe test for the KKT-Limiter and semi-smooth Newton algorithm. The solution $u_h$ and values of the Lagrange multiplier are shown in Figure \ref{Buckley_Leverett3} and the derivative $Q_h$ in Figure \ref{Buckley_Leverett4}. The results on the two meshes compare well and the limiter ensures that the positivity and maximum bounds are satisfied.

The number of Newton iterations necessary to obtain a minimum value $10^{-8}$ for the residual $F(z)$ and Newton update $d$ in Algorithm \ref{algorithm_1}  to stop the Newton iterations for each DIRK stage strongly varies. It depends on the type of equation, time-step and nonlinearity. In general, the time step is chosen such that the number of Newton iterations for each DIRK stage is between 5 and 20. For most time dependent problems the CFL number is then close to one, which is necessary to ensure time-accuracy. Only for the Buckley-Leverett equation with gravity the time step frequently had to be less than one in order to deal with the strong nonlinearity of the problem. In the computations we did not observe a minimum time step to ensure positivity as noticed in \cite{qin2018implicit}.

\section{Conclusions}
\label{sec:conclusions}
In this paper we present a novel framework to combine positivity preserving limiters for discontinuous Galerkin discretizations with implicit time integration methods. This approach does not depend on the specific type of discontinuous Galerkin discretization and is also applicable to e.g. finite volume discretizations. The key features of the numerical method is the formulation of the positivity constraints as a Karush-Kuhn-Tucker problem and the development of an active set semi-smooth Newton method that accounts for the non-smoothness of the algebraic equations. The algorithm was successfully tested on a number of increasingly difficult test cases, which required that the positivity constraints are satisfied in order to obtain meaningful results. The KKT Limiter does not negatively affect the accuracy for smooth problems and accurately preserves the positivity constraints. Future work will focus on the extension of the KKT Limiter to ensure also monotonicity of the solution.

\appendix
\section{Derivation of Clarke directional derivative}\label{Appendix1}
For completeness we give here a derivation of the terms  \cref{Clarke_direc_derivb,Clarke_direc_derivc} in the Clarke directional derivative of $F(z)$ in \cref{KKTnonlineareq}. We will follow the approach outlined in \cite{ito2008lagrange}. Define $z:=(x,\mu,\lambda), \bar{z}:=(\bar{x},\bar{\mu},\bar{\lambda}), d:=(u,v,w)\in\mathbb{R}^{p}$, with $p=n+l+m$. Consider  $\bar{F}(z)=F_{i+n+l}(z)$,  $i\in\beta(z)$. The other Clarke directional derivatives of $F$ are straightforward to compute.  If we consider \cref{thetaclark} only for the  contribution of $\bar{F}(z)$ to the merit function to $\theta(z)$ and use \cref{KKTnonlineareq} and a Taylor expansion of $\bar{F}(z)$ around $z$, then we obtain
\begin{align*}
\bar{\theta}^0(z;d)&=\limsup_{\bar{z}\rightarrow z, t\downarrow 0^+}\frac{1}{t}\Big( \bar{F}(z),\min(-g(\bar{x}+tu),\bar{\lambda}+tw)-\min(-g(\bar{x}),\bar{\lambda})\Big)\\
&=\limsup_{\bar{z}\rightarrow z, t\downarrow 0^+}\frac{1}{t}\Big( \bar{F}(z),\min(-g(x)-J(\bar{x}+tu-x),\bar{\lambda}+tw)\\
&\hspace*{77pt}-\min\big(-g(x)-J(\bar{x}-x),\bar{\lambda}\big)\Big),
\end{align*}
with $J:=D_xg(x)\in\mathbb{R}^{m\times n}$.
Here, higher order terms are omitted since they will become zero in the limit. Define $h(x):=-g(x)+Jx$, then
\begin{align}
\bar{\theta}^0(z;d)&=\limsup_{\bar{z}\rightarrow z, t\downarrow 0^+}\frac{1}{t}\Big( \bar{F}(z),\min(-J\bar{x}-tJu+h(x),\bar{\lambda}+tw)\label{theta0deriv}\\
&\hspace*{77pt}-\min(-J\bar{x}+h(x),\bar{\lambda})\Big).\nonumber
\end{align}
For $u\in\mathbb{R}^n, w\in\mathbb{R}^m$, define $r\in\mathbb{R}^m$ by
\begin{subequations}
\begin{alignat}{2}
&r_i<0&\;\text{on}\;S_1:=&\{i\in\beta(z)\;\vert\;\bar{F}_i(z)>0,-(Ju)_i>w_i\}\nonumber\\
&&\cup&\{i\in\beta(z)\;\vert\;\bar{F}_i(z)\leq 0,-(Ju)_i\leq w_i\},\label{r-conditionsA}\\[5pt]
&r_i>0&\;\text{on}\;S_2:=&\{i\in\beta(z)\;\vert\;\bar{F}_i(z)>0,-(Ju)_i\leq w_i\}\nonumber\\
&&\cup&\{i\in\beta(z)\;\vert\;\bar{F}_i(z)\leq 0,-(Ju)_i> w_i\}.\label{r-conditionsB}
\end{alignat}\label{r-conditions}
\end{subequations}
Let $\bar{x}\in\mathbb{R}^n$ be such that
\begin{equation}
-J\bar{x}+h(x)=\bar{\lambda}+r.\label{xbar_equation}
\end{equation}
Note, such an $\bar{x}$ exists for $i\in\beta(z)$ since \cref{xbar_equation} is equivalent with $-Ju=w+r$ with $u=\bar{x}-x$ and $w=\bar{\lambda}-\lambda$ as components of the search direction $d$.
Choose $t\in(0,t_{\bar{x}})$  for $t_{\bar{x}}>0$ such that
\begin{subequations}
\begin{alignat}{2}
(-J\bar{x}+h(x)-tJu)_i&<(\bar{\lambda}+tw)_i&\quad&\text{for}\;i\in S_1,\\
(-J\bar{x}+h(x)-tJu)_i&>(\bar{\lambda}+tw)_i&\quad&\text{for}\;i\in S_2.
\end{alignat}\label{S-conditions}
\end{subequations}
Note, such a $t_{\bar{x}}$ exists, see Remark \ref{remark1}.
We then obtain
\begin{equation*}
\min((-J\bar{x}+h(x)-tJu)_i,(\bar{\lambda}+tw)_i)=
\begin{cases}(-J\bar{x}+h(x)-tJu)_i\quad&\text{for}\;i\in S_1,\\
(\bar{\lambda}+tw)_i\quad&\text{for}\;i\in S_2.
\end{cases}
\end{equation*}
Use now \cref{xbar_equation} and \cref{r-conditions} then
\begin{align*}
\min((-J\bar{x}+h(x))_i,\bar{\lambda}_i)&=\min(\bar{\lambda}_i+r_i,\bar{\lambda}_i)=\begin{cases}
\bar{\lambda}_i+r_i&\text{for}\; i\in S_i,\\
\bar{\lambda}_i&\text{for}\; i\in S_2.
\end{cases}
\end{align*}
Combining the above results and using \cref{xbar_equation} again gives
\begin{align*}
\min((-J\bar{x}+h(x)-tJu)_i,(\bar{\lambda}+tw)_i)&-\min((-J\bar{x}+h(x))_i,\bar{\lambda}_i)\\[5pt]
&=\begin{cases}
-t(Ju)_i\quad&\text{for}\;i\in S_1,\\
tw_i,&\text{for}\; i\in S_2,
\end{cases}\\[5pt]
&=\begin{cases}
t\max(-(Ju)_i,w_i)\quad&\text{if}\;\bar{F}_i(z)>0,\\
t\min(-(Ju)_i,w_i),&\text{if}\;\bar{F}_i(z)\leq0.
\end{cases}
\end{align*}
Taking the limit in \cref{theta0deriv} and using \cref{thetader1} for  $\bar{\theta}(z;d)$ then gives \cref{Clarke_direc_derivb,Clarke_direc_derivc}.

\begin{remark}\label{remark1}
Conditions \cref{r-conditions} imply \cref{S-conditions}. Use $-J\bar{x}+h(x)=\bar{\lambda}+r$ in \cref{S-conditions}, then we obtain
\begin{alignat}{2}
(r-tJu)_i&<tw_i&\quad&\text{for}\; i\in S_1,\label{con1}\\
(r-tJu)_i&>tw_i&\quad&\text{for}\; i\in S_2\label{con2}.
\end{alignat}

I. If $i\in S_1$, $\bar{F}_i(z)>0$ then from \cref{r-conditionsA} we obtain $-(Ju)_i-w_i>0$ and  \cref{con1} implies $r_i+t(-(Ju)_i-w_i)<0$. Choose $t<\frac{-r_i}{-(Ju)_i-w_i}=t_{\bar{x}}$. Since $r_i<0$ and $-(Ju)_i-w_i>0$ for $i\in S_1$, $\bar{F}_i(z)>0$ we obtain that $t_{\bar{x}}>0$.

II. If $i\in S_1$, $\bar{F}_i(z)\leq 0$, then \cref{r-conditionsA} implies $-(Ju)_i-w_i\leq 0$ and  \cref{con1} gives $r_i+t(-(Ju)_i-w_i)<0$. Since both $r_i$ and $-(Ju)_i-w_i<0$ and any $t>0$ will imply \cref{con1}.

The proof for $i\in S_2$ is completely analogous are therefore omitted. Hence there exists a $t_{\bar{x}}>0$ for \cref{S-conditions}.
\end{remark}
\section{Verification of conditions for quasi-directional derivative}\label{Appendix2}
In this section we show that the quasi-directional derivative \cref{quasidirectionDeriv1} satisfies the conditions stated in \cref{quasidirectional_derivative}, which are necessary to ensure converge of the Newton algorithm defined in Algorithm \ref{Newton_algorithm}.

Consider condition \cref{quasi_der_cond1}: First note that
\begin{alignat*}{2}
&F_i^\prime(z;d)=F_i^0(z;d)=G_i(z;d),&\quad &i\in N_n,\\
&F_{i+n}^\prime(z;d)=F_{i+n}^0(z;d)=G_{i+n}(z;d),&\quad &i\in N_l,\\
&F_{i+n+l}^\prime(z;d)=F_{i+n+l}^0(z;d)=G_{i+n+l}(z;d),&\quad &i\in\alpha_\delta(z)\cup\gamma_\delta(z),
\end{alignat*}
since $\alpha_\delta(z)\cup\gamma_\delta(z)\subset\alpha(z)\cup\gamma(z)$. If $i\in\beta_\delta(z)$ and $F_{i+n+l}(z)\leq 0$ then
\begin{equation*}
\min(-(Ju)_i,w_i)\leq-(Ju)_i,w_i.
\end{equation*}
Since $F_{i+n+l}(z)\leq 0$ this implies
\begin{equation*}
F_{i+n+l}(z)\min(-(Ju)_i,w_i)\geq F_{i+n+l}(z)(-(Ju)_i),F_{i+n+l}(z)w_i.
\end{equation*}
If $i\in\beta_\delta(x)$ and $F_{i+n+l}(z)>0$ then
\begin{equation*}
-(Ju)_i,w_i\leq \max(-(Ju)_i,w_i).
\end{equation*}
Hence, since $F_{i+n+l}(z)>0$ this implies
\begin{equation*}
F_{i+n+l}(z)(-(Ju)_i), F_{i+n+l}(z)w_i\leq F_{i+n+l}(z)\max(-(Ju)_i,w_i).
\end{equation*}
 Comparing all terms then immediately shows that $G(z;d)$ satisfies \cref{quasi_der_cond1} and \cref{Gupperbound}. Condition \cref{quasi_der_cond2}  directly follows from the definition of $G$ in \cref{quasidirectional_derivative}.

\section*{Acknowledgments}
We would like to acknowledge Mrs. Fengna Yan from USTC and the University of Twente for her contributions in testing the KKT-Limiter for several DG discretizations.

\bibliographystyle{siamplain}
\bibliography{JJWvdVegt_bibreferences}
\end{document}